# A SURVEY OF MAX-TYPE RECURSIVE DISTRIBUTIONAL EQUATIONS

By David J. Aldous[1] and Antar Bandyopadhyay

*University of California, Berkeley and University of Minnesota*

In certain problems in a variety of applied probability settings (from probabilistic analysis of algorithms to statistical physics), the central requirement is to solve a *recursive distributional equation* of the form $X \stackrel{d}{=} g((\xi_i, X_i), i \geq 1)$. Here $(\xi_i)$ and $g(\cdot)$ are given and the $X_i$ are independent copies of the unknown distribution $X$. We survey this area, emphasizing examples where the function $g(\cdot)$ is essentially a "maximum" or "minimum" function. We draw attention to the theoretical question of *endogeny*: in the associated recursive tree process $X_{\mathbf{i}}$, are the $X_{\mathbf{i}}$ measurable functions of the innovations process $(\xi_{\mathbf{i}})$?

**1. Introduction.** Write $\mathcal{P}$ for the space of probability distributions on a space $S$; in our examples, $S$ will usually be $\mathbb{R}$ or a subset of $\mathbb{R}$. Suppose we are given a joint distribution for some family of random variables $(\xi_i, i \geq 1)$, and given an $S$-valued function $g(\cdot)$ with appropriate domain (Section 2.1 gives this setup more carefully). Then we can define a map $T : \mathcal{P} \to \mathcal{P}$ as follows:

> $T(\mu)$ is the distribution of $g((\xi_i, X_i), i \geq 1)$, where the $(X_i)$ are independent with distribution $\mu$, independent of the family $(\xi_i)$.

Within this general framework one can ask about existence and uniqueness of *fixed points*, that is, distributions $\mu$ such that $T(\mu) = \mu$, and about *domain of attraction* for a fixed point $\mu$, that is, for what initial $\nu$ do we have $T^n(\nu) \to \mu$ as $n \to \infty$. One can rewrite such a fixed point equation in terms of random variables as

$$(1) \qquad X \stackrel{d}{=} g((\xi_i, X_i), i \geq 1)$$

Received February 2004; revised June 2004.

[1]Supported by NSF Grant DMS-02-03062.

*AMS 2000 subject classifications.* 60E05, 62E10, 68Q25, 82B44.

*Key words and phrases.* Branching process, branching random walk, cavity method, coupling from the past, fixed point equation, frozen percolation, mean-field model of distance, metric contraction, probabilistic analysis of algorithms, probability distribution, probability on trees, random matching.

where the independence property is assumed implicitly. We introduce the phrase *recursive distributional equation* (RDE) for equations of format (1), as opposed to alternate kinds of fixed point equation. RDEs have arisen in a variety of settings:

- Galton–Watson branching processes and related random trees,
- probabilistic analysis of algorithms with suitable recursive structure,
- statistical physics models on trees,
- statistical physics and algorithmic questions in the mean-field model of distance.

Three aspects of this topic have been well studied. Perhaps the best known fixed point equation is

$$(2) \qquad\qquad X \stackrel{d}{=} 2^{-1/2}(X_1 + X_2) \qquad (S = \mathbb{R})$$

whose solutions are the Normal$(0, \sigma^2)$ family. This example extends to give characterizations of stable distributions [63]. Moreover, there is a classical topic "characterization of probability distributions" [43] which considers the named families of distributions in mathematical statistics and studies many different types of characterization, some of which can be put into form (1). But this aspect is rather remote from our concerns. A second aspect concerns general methods for establishing existence or uniqueness of fixed points. Two natural methods (monotonicity; metric contraction) are recalled in Section 2.2, while the more elaborate method of "a.s. unique tree representations" or "tree-structured coupling from the past" is described in Section 2.6. The third aspect is the *linear* subcase $g((\xi_i, X_i), i \geq 1) = \sum_i \xi_i X_i$ and its variants, which we review in Section 3. This is well understood for $S = \mathbb{R}^+$, though not so well understood for $S = \mathbb{R}$.

The main purpose of this survey article is to illustrate the variety of contexts where RDEs have arisen, and to draw attention to another subclass of RDEs, those involving *max-type* functions $g$. We present in Sections 4–7 a collection of around ten examples (summarized in Table 1) of max-type RDEs arising from concrete questions. Most have been studied in detail elsewhere; in each case we seek to outline the underlying problem, describe how it leads to an RDE, and give information about solutions of general or special cases. Section 8 provides brief remarks on Monte Carlo methods, process analogs and continuous-time analogs, and lists the numerous open problems and conjectures.

On the theoretical side, in Section 2.3 we set out carefully some general theory, rather natural from the statistical physics or interacting particle system viewpoint but less apparent from the algorithms viewpoint, which relates RDEs to a type of tree-indexed process $(X_\mathbf{i})$ which we call *recursive tree processes* (RTPs). In particular we introduce the *endogenous* property (Definition 7), that in an RTP $X_\mathbf{i}$ is a measurable function of the driving tree-indexed process $(\xi_\mathbf{i})$ without any external randomness being needed, and show (Theorem 11) that endogeny is equivalent to a *bivariate uniqueness* property.

TABLE 1
*Some max-type RDEs. Functions $g(\cdot)$ for which the RDE $X \overset{d}{=} g((\xi_i, X_i), i \geq 1)$ are discussed*[*]

| Section | $g(\cdot)$ | Underlying model | Endog? | Comments |
|---|---|---|---|---|
| | $S = \mathbb{R}^+$ | | | |
| 4.2 | $\max_i(X_i + \xi_i)^+$ | Range of BRW | Yes | |
| 4.3 | $\min_i(X_i + \xi_i)^+$ | Algorithm for BRW range | Yes | |
| 4.6 | $\max_i(\xi_i - X_i)^+$ | Matching on GW tree | Yes | |
| 4.4 | $\xi_0 + \max_i(\xi_i X_i)$ | Discounted tree sums | Yes | $\xi_0 = 0$ reduces to BRW extremes |
| 4.4 | $\xi_0 + \min_i(\xi_i X_i)$ | Discounted tree sums | Yes | See (49) |
| 4.6 | $(\xi_0 - \sum_i X_i)^+$ | Independent subset GW tree | Yes | |
| 7.2 | $\sum_i(c - \xi_i + X_i)^+$ | Percolation of MSTs | Yes | Determines critical $c$ |
| 7.6 | See (98) | First passage percolation | Conj. $Y$ | Mean-field scaling analysis |
| | $S = \mathbb{R}$ | | | |
| 5 | $c + \max_i(X_i + \xi_i)$ | Extremes in BRW | No | $c$ specified by dist($\xi_i$) |
| 7.3 | $\min_i(\xi_i - X_i)$ | Mean-field minimal matching | Yes | |
| 7.4 | $\min^{[2]}_i(\xi_i - X_i)$ | Mean-field TSP | Conj. $Y$ | $\min^{[2]}$ denotes second smallest |
| | Other $S$ | | | |
| 6 | $\Phi(\min(X_1, X_2), \xi_0)$ | Frozen percolation on tree | Yes | $\Phi$ defined in Section 6 |
| 7.6 | See (96), (97), (98) | Mean-field scaling | Conj. $Y$ | $S = \mathbb{R}^2$ or $\mathbb{R}^3$ |

[*]Note $x^+ = \max(x, 0)$. For $S = \mathbb{R}$ a "max" problem is equivalent to a "min" problem by transforming $X$ to $-X$, but for $S = \mathbb{R}^+$ this does not work: the problems in Sections 4.2 and 4.3 are different. Typically the $(\xi_i)$ are either i.i.d. or are the successive points of a Poisson process on $(0, \infty)$. "Endogenous" refers to fundamental solution. Key to acronyms: BRW, branching random walk; GW, Galton–Watson; MST, minimal spanning tree; TSP, traveling salesman problem.

A concluding Section 9 will attempt to review the big picture.

1.1. *Three uses of RDEs.* When we look at how RDEs arise within specific models in the Table 1 examples, we will see three broad categories of use, which seem worth mentioning at the start.

1.1.1. *Direct use.* Here is the prototype example of *direct use*, where the original question asks about a random variable $X$ and the distribution of $X$ itself satisfies an RDE.

EXAMPLE 1. Let $X$ be the total population in a Galton–Watson branching process where the number of offspring is distributed as $\xi$. In the case $\mathbb{E}\xi \leq 1$ [and $\mathbb{P}(\xi = 1) \neq 1$] it is well known that $X < \infty$ a.s., and then easy to check that $\mathrm{dist}(X)$ is the unique solution of the RDE

$$X \stackrel{d}{=} 1 + \sum_{i=1}^{\xi} X_i \qquad (S = \mathbb{Z}^+).$$

We will see other direct uses in Proposition 25 and in the examples in Section 4.4.

1.1.2. *Indirect use.* The simplest kind of indirect use is where the quantity of interest can be written in terms of known quantities and some other quantity which can be analyzed via an RDE. See Proposition 28 and Theorem 41 for results of this kind. But there is a more intriguing kind of indirect use which we call a $540°$ *argument*, exemplified in the frozen percolation model of Section 6 and also used in the mean-field combinatorial optimization problems in Sections 7.3–7.6. In these examples we start with a heuristically defined quantity $X$, and a heuristic argument that it should satisfy an RDE. Next we make a rigorous argument by first solving the RDE and then using the associated, rigorously defined RTP as building blocks for a rigorous construction.

1.1.3. *Critical points and scaling laws.* We introduce this idea with an artificial example.

EXAMPLE 2. Let $\xi$ be $\mathbb{R}$-valued, $\mathbb{E}\xi = \beta$, and let $(\xi_i, i \geq 1)$ be independent copies of $\xi$. For fixed $c \in \mathbb{R}$ consider the RDE

(3)                     $$X \stackrel{d}{=} \max(0, X + \xi - c) \qquad (S = \mathbb{R}^+).$$

Then there is a solution $X_c$ on $\mathbb{R}^+$ if and only if $c > \beta$. Moreover, if $\mathrm{var}(\xi) \in (0, \infty)$, then

(4)                     $$\mathbb{E}X_c \sim \frac{\mathrm{var}(\xi)}{2(c - \beta)} \qquad \text{as } c \downarrow \beta.$$

Here (3) is a *Lindley equation* from classical queuing theory [12], and it is straightforward that for $c > \beta$ the solution is

(5)                     $$X_c \stackrel{d}{=} \max_{j \geq 0} \sum_{i=1}^{j} (\xi_i - c).$$

This $X_c$ is a.s. finite by the *strong law of large numbers*, and the stated asymptotics (4) follow from, for example, weak convergence of random walks to Brownian motion with drift.

We will see later three examples of problems involving *critical values* or *near-critical behavior* of some random system. In such problems there is a parameter $c$ and we are interested in a critical value $c_{\mathrm{crit}}$ of $c$ defined as the value where some "phase transition" occurs, or in behavior of the system for $c$ near $c_{\mathrm{crit}}$. In Section 7.2 we see an example where the critical point is determined as the boundary between existence and nonexistence of a solution to an RDE (Proposition 56). In Section 4.3 we see how aspects of near-critical behavior may be reduced to study of near-critical solutions of an RDE (Theorem 29), and Section 7.6 contains more sophisticated variations on that theme. Note that in Example 2, result (4) shows that the behavior of solutions near the critical point scales in a simple way that does not depend on the details of the distribution of $\xi$; according to the statistical physics paradigm of *universality* one should expect such scaling laws to arise in most natural problems.

1.2. *The cavity method.* One particular topic of current interest concerns the *cavity method* in statistical physics, applied in the context of combinatorial optimization in mean-field settings. There is a methodology for seeking rigorous proofs, in which the central issue becomes uniqueness of solution of some problem-dependent RDE. We will elaborate slightly in Section 7.5.

## 2. The general setting.

2.1. *A precise setup.* Here we record a more careful setup for RDEs. Let $(S, \mathcal{S})$ be a measurable space, and let $\mathcal{P}(S)$ be the set of probability measures on $(S, \mathcal{S})$. Let $(\Theta, \mathcal{T})$ be another measurable space. Construct

$$\Theta^* := \Theta \times \bigcup_{0 \le m \le \infty} S^m,$$

where the union is a disjoint union and where $S^m$ is product space, interpreting $S^\infty$ as the usual infinite product space and $S^0$ as a singleton set, which we will write as $\{\Delta\}$. Let $g : \Theta^* \to S$ be measurable. Let $\nu$ be a probability measure on $\Theta \times \bar{\mathbb{Z}}^+$, where $\bar{\mathbb{Z}}^+ := \{0, 1, 2, \ldots; \infty\}$. These objects can now be used to define a measurable map $T : \mathcal{P}(S) \to \mathcal{P}(S)$ as follows. Write $\le^* N$ to mean $\le N$ for $N < \infty$ and to mean $< \infty$ for $N = \infty$.

DEFINITION 3. $T(\mu)$ is the distribution of $g(\xi, X_i, 1 \le i \le^* N)$, where:

(i) $(X_i, i \ge 1)$ are independent with distribution $\mu$;
(ii) $(\xi, N)$ has distribution $\nu$;
(iii) the families in (i) and (ii) are independent.

Equation (1) fits this setting by writing $\xi = (\xi_i)$. In most examples there is a sequence $(\xi_i)$, but for theoretical discussion we regard such a sequence as a single random element $\xi$.

In examples where $P(N = \infty) > 0$ a complication often arises. It may be that $g(\cdot)$ is not well defined on all of $\Theta \times S^\infty$, although $g(\xi, X_i, 1 \leq i \leq^* N)$ is well defined almost surely for $(X_i)_{i \geq 1}$ i.i.d. with distribution in a restricted class of probabilities on $S$. For such examples and also for other cases where it is natural to restrict attention to distributions satisfying some conditions (like moment conditions), we allow the more general setting where we are given a subset $\mathcal{P} \subseteq \mathcal{P}(S)$ such that $g(\xi, X_i, 1 \leq i \leq^* N)$ is well defined almost surely for i.i.d. $(X_i)_{i \geq 1}$ with distribution in $\mathcal{P}$. Now $T$ is well defined as a map

$$(6) \qquad\qquad T \colon \mathcal{P} \to \mathcal{P}(S).$$

In this *extended* case it is natural to seek, but maybe hard to find, a subset $\mathcal{P}$ such that $T$ maps $\mathcal{P}$ *into* $\mathcal{P}$.

2.2. *Monotonicity and contraction.* There are standard tools for studying maps $T \colon \mathcal{P}(S) \to \mathcal{P}(S)$ which do not depend on the map arising in the particular way of Definition 3. First suppose $S \subset \bar{\mathbb{R}}$ is an interval of the form $[0, x_0]$ for some $x_0 < \infty$, or $S = [0, \infty)$. Consider the usual *stochastic partial order* $\preceq$ on $\mathcal{P}(S)$:

$$\mu_1 \preceq \mu_2 \quad \text{iff } \mu_1[0, x] \geq \mu_2[0, x], \qquad x \in S.$$

Say $T$ is *monotone* if

$$\mu_1 \preceq \mu_2 \quad \text{implies} \quad T(\mu_1) \preceq T(\mu_2).$$

Note that, writing $\delta_0$ for the probability measure degenerate at 0, if $T$ is monotone, then the sequence of iterates $T^n \delta_0$ is increasing, and then the limit

$$\lim_n T^n \delta_0 = \mu_*$$

exists in the sense of weak convergence on the compactified interval $[0, \infty]$.

LEMMA 4 (Monotonicity lemma). *Let $S$ be an interval as above. Suppose $T$ is monotone. If $\mu_*$ gives nonzero measure to $\{\infty\}$, then $T$ has no fixed point on $\mathcal{P}(S)$. If $\mu_*$ gives zero measure to $\{\infty\}$, and if $T$ is continuous with respect to increasing limits $[\mu_n \uparrow \mu_\infty$ implies $T(\mu_n) \uparrow T(\mu_\infty)]$, then $\mu_*$ is a fixed point of $T$, and $\mu_* \preceq \mu$, for any other fixed point $\mu$.*

This obvious result parallels the notion of *lower invariant measure* in interacting particle systems [47].

Returning to the case of general $S$, the Banach contraction theorem specializes to LEMMA 5 (The contraction method). *Let $\mathcal{P}$ be a subset of $\mathcal{P}(S)$ such that $T$ maps $\mathcal{P}$ into $\mathcal{P}$. Let $d$ be a complete metric on $\mathcal{P}$. Suppose $T$ is a (strict) contraction, that is,*

$$\sup_{\mu_1 \neq \mu_2 \in \mathcal{P}} \frac{d(T(\mu_1), T(\mu_2))}{d(\mu_1, \mu_2)} < 1.$$

*Then $T$ has a unique fixed point $\mu$ in $\mathcal{P}$, whose domain of attraction is all of $\mathcal{P}$.*

A thorough account of specific metrics can be found in [57]. Most commonly used is the Wasserstein metric on distributions on $\mathbb{R}$ with finite $p$th moment, $1 \leq p < \infty$:

$$(7) \qquad d_p(\mu, \nu) := \inf\{(\mathbb{E}[|Z - W|^p])^{1/p} | Z \stackrel{d}{=} \mu \text{ and } W \stackrel{d}{=} \nu\}.$$

Contraction is a powerful tool in the "linear" case of Section 3, where it also provides rates of convergence in the context of probabilistic analysis of algorithms. For max-type operations it seems less widely useful (see, e.g., the remark below Open Problem 62) except in simple settings (e.g., Theorem 32). It is also worth mentioning that in several examples we have no rigorous proofs of existence of fixed points (see Sections 5 and 7.4) and so the use of other fixed point theorems from analysis [45] might be worth exploring.

2.3. *Recursive tree processes.* Consider again the setup from Section 2.1. Rather than considering only the induced map $T$, one can make a richer structure by interpreting

$$X = g(\xi, X_i, 1 \leq i \leq^* N)$$

as a relationship between random variables. In brief, we regard $X$ as a value associated with a "parent" which is determined by the values $X_i$ at $N$ "children" and by some "random noise" $\xi$ associated with the parent. One can then extend to grandchildren, great grandchildren and so on in the obvious way. We write out the details carefully in the rest of this section.

Write $\mathbb{T}$ for the set of all possible descendants $\mathbf{i}$, where $\mathbf{i} = i_1 i_2 \cdots i_d$ denotes a $d$th-generation individual, the $i_d$th child of its parent $i_1 i_2 \cdots i_{d-1}$. Label the root as $\varnothing$. Make $\mathbb{T}$ a tree by adding parent–child edges. Given the distribution $\nu$ on $\Theta \times \bar{\mathbb{Z}}$ from Section 2.1, for each $\mathbf{i} \in \mathbb{T}$ let $(\xi_{\mathbf{i}}, N_{\mathbf{i}})$ have distribution $\nu$, independently as $\mathbf{i}$ varies. Recall the function $g$ from Section 2.1. This structure—the random pairs $(\xi_{\mathbf{i}}, N_{\mathbf{i}})$, $\mathbf{i} \in \mathbb{T}$, which are i.i.d. $(\nu)$, and the function $g$—we call a *recursive tree framework* (RTF). In the setting of an RTF suppose that, jointly with the random objects above, we can construct $S$-valued random variables $X_{\mathbf{i}}$ such that for each $\mathbf{i}$

$$(8) \qquad X_{\mathbf{i}} = g(\xi_{\mathbf{i}}, X_{\mathbf{i}j}, 1 \leq j \leq^* N_{\mathbf{i}}) \qquad \text{a.s.}$$

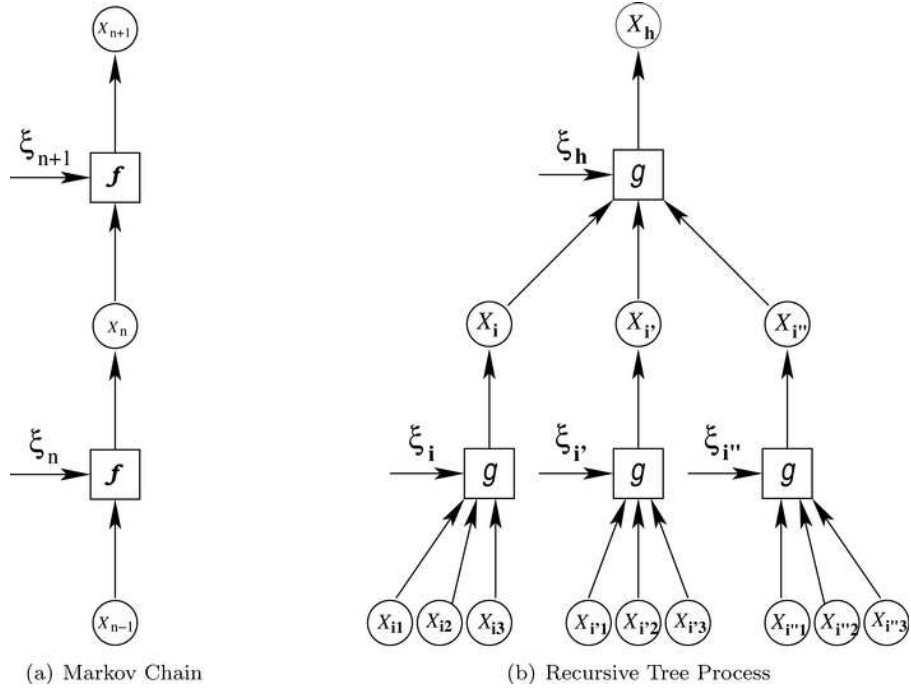

(a) Markov Chain                          (b) Recursive Tree Process

Fig. 1. *On the left is a Markov chain represented as an iterated function system: $X_n$ is the "output" of $f$ with "inputs" $\xi_n$ and $X_{n-1}$. On the right is an RTP; $X_{\mathbf{i}}$ is the "output" of $g$ with inputs $\xi_{\mathbf{i}}$ and $(X_{\mathbf{j}}, \mathbf{j}$ child of $\mathbf{i})$. In the figure, $\mathbf{h}$ is the parent of $\mathbf{i}$ and $\mathbf{i}', \mathbf{i}'', \ldots$ are siblings of $\mathbf{i}$.*

and such that, independent of the values of $\{\xi_{\mathbf{i}}, N_{\mathbf{i}} | \mathbf{i}$ in first $d-1$ generations$\}$, the random variables $\{X_{\mathbf{i}} | \mathbf{i}$ in generation $d\}$ are i.i.d. with some distribution $\mu_d$. Call this structure (an RTF jointly with the $X_{\mathbf{i}}$) a *recursive tree process* (RTP). If the random variables $X_{\mathbf{i}}$ are defined only for vertices $\mathbf{i}$ of depth $\leq d'$, then call it an RTP of depth $d'$. See Figure 1.

Now an RTF has an *induced* map $T \colon \mathcal{P}(S) \to \mathcal{P}(S)$ as in Definition 3. [In the extended case (6) we need to assume that $T$ maps $\mathcal{P}$ into $\mathcal{P}$.] Note that the relationship between an RTF and an RTP mirrors the relationship between a Markov transition kernel and a Markov chain. Fix an RTF. Given $d$ and an arbitrary distribution $\mu^0$ on $S$, there is an RTP of depth $d$ in which the generation-$d$ vertices are defined to have distribution $\mu_d = \mu^0$. Then the distributions $\mu_d, \mu_{d-1}, \mu_{d-2}, \ldots, \mu_0$ at decreasing generations $d, d-1, d-2, \ldots, 0$ of the tree are just the successive iterates $\mu^0, T(\mu^0), T^2(\mu^0), \ldots, T^d(\mu^0)$ of the map $T$. Figures 1 and 2 attempt to show the analogy between RTPs and Markov chains.

One should take a moment to distinguish RTPs from other structures involving tree-indexed random variables. For instance, a *branching Markov*

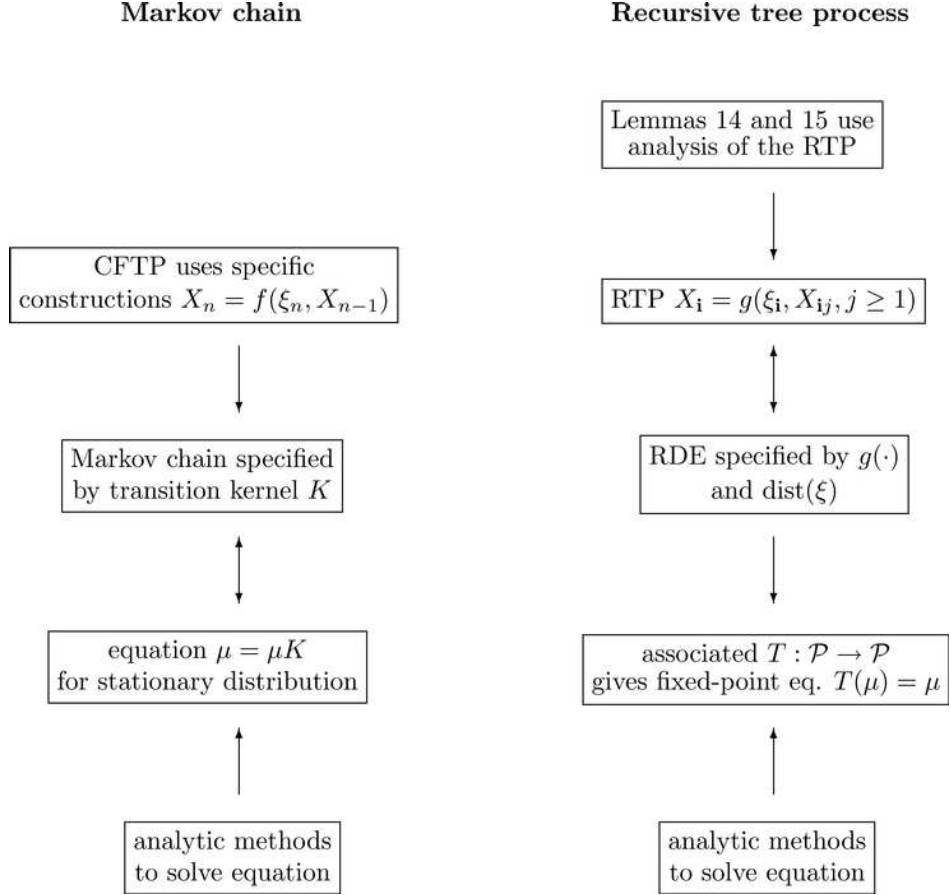

**Markov chain**                                    **Recursive tree process**

Fig. 2.   *The analogy between Markov chains and RTPs. Specifying a transition kernel is exactly what is needed to write the equation for a stationary distribution, and exactly what is needed to specify the chain. Analogously, specifying $S$ and $g(\cdot)$ and dist($\xi$) is exactly what is needed to write the RDE, and exactly what is needed to specify the RTP. But note these equivalences ↕ occur at different conceptual levels.*

*chain* can also be represented as a family $(X_{\mathbf{i}})$. But its essential property is that, conditional on the value $X_{\mathbf{i}}$ at a parent $\mathbf{i}$, the values $(X_{\mathbf{i}1}, X_{\mathbf{i}2}, \ldots)$ at the children $\mathbf{i}1$, $\mathbf{i}2, \ldots$ are i.i.d. An RTP in general does not have this property. Conceptually, in branching processes one thinks of the "arrow of time" as pointing away from the root, whereas in an RTF the arrow points toward the root.

Call an RTP *invariant* if the marginal distributions of $X_{\mathbf{i}}$ are identical at all depths. We have the following obvious analog of Markov chain stationarity.

LEMMA 6.   *Consider an RTF. A distribution $\mu$ is a fixed point of the induced map $T$ if and only if there is an invariant RTP with marginal distributions $\mu$.*

An invariant RTF could be regarded as a particular case of a *Markov random field*, but the special "directed tree" structure of RTFs makes them worth distinguishing from general Markov random fields.

A central theme of this survey paper is that for certain purposes, the most useful way of interpreting an RDE (1) is as the defining identity for an invariant RTP.

2.4. *Endogeny and bivariate uniqueness.*   Now imagine (8) as a system of equations for "unknowns" $X_{\mathbf{i}}$ in terms of "known data" $\xi_{\mathbf{i}}$. It is natural to ask if the solution depends only on the data. We formalize this as the following *endogenous* property. Write

$$(9) \qquad\qquad \mathcal{G}_{\mathbb{T}} = \sigma(\xi_{\mathbf{i}}, N_{\mathbf{i}}, \ \mathbf{i} \in \mathbb{T}).$$

DEFINITION 7.   An invariant RTP is called *endogenous* if

$$X_{\varnothing} \text{ is } \mathcal{G}_{\mathbb{T}}\text{-measurable.}$$

A rephrasing is more intuitive. Within an RTF there is an embedded Galton–Watson tree $\mathcal{T}$ rooted at $\varnothing$, whose offspring distribution $N$ is just the marginal in $\nu = \mathrm{dist}(\xi, N)$. That is, the root $\varnothing$ has $N_{\varnothing}$ children; each such child $\mathbf{i}$ has $N_{\mathbf{i}}$ children, and so on; $\mathcal{T}$ is the random set of all such descendants of the root $\varnothing$. Write

$$(10) \qquad\qquad \mathcal{G} = \sigma(\xi_{\mathbf{i}}, N_{\mathbf{i}}, \ \mathbf{i} \in \mathcal{T}).$$

Then endogeny is equivalent to

$$X_{\varnothing} \text{ is } \mathcal{G}\text{-measurable}$$

and this is the criterion we use in practice.

It is intuitively clear (and true: Lemma 14) that when the Galton–Watson tree $\mathcal{T}$ is a.s. finite there will be a unique invariant RTP and it will be endogenous. But when $\mathcal{T}$ is infinite the "boundary behavior" may cause uniqueness and/or endogeny to fail. Theorem 11 will show that the endogenous property is equivalent to a certain *bivariate uniqueness* property. The simple artificial Examples 8, 10 and 13 should help to distinguish these properties from the "unique fixed point of $T$" property.

Our first example shows that one cannot tell whether or not the endogenous property holds just by looking at $T$, even when the fixed point is unique. Write Bern($p$) for the Bernoulli($p$) distribution on $\{0, 1\}$.

Example 8.   Take $S = \{0,1\}$. Define $T : \mathcal{P}(S) \to \mathcal{P}(S)$ by $T(\mu) = \text{Bern}(1/2)$ for all $\mu$. So $\text{Bern}(1/2)$ is the unique fixed point. We will give two RTPs with this induced $T$, one satisfying and the other failing the endogenous property.

First take $(\xi, N)$ with $N = 2$ and $\xi \overset{d}{=} \text{Bern}(1/2)$, and $g(a, x_1, x_2) = a$. Clearly the induced $T$ is as above. In the associated RTP where $X_{\mathbf{i}}$ has $\text{Bern}(1/2)$ marginals, observe that $X_{\varnothing} = \xi_{\varnothing}$ and so the endogenous property holds. Now consider the well-known von Neumann random bit extractor [30], which is a certain function $\bar{g} : \{0,1\}^{\infty} \to \{0,1\}$ which, applied to an independent $\text{Bern}(p)$ input sequence for any $0 < p < 1$, gives a $\text{Bern}(1/2)$ output. Set

$$g(a, x_1, x_2, \dots) = \begin{cases} a, & \text{if } x_1 = x_2 = x_3 \cdots, \\ \bar{g}(x_1, x_2, \dots), & \text{if not.} \end{cases}$$

Take $(\xi, N)$ with $N = \infty$ and $\xi \overset{d}{=} \text{Bern}(1/2)$, and then the induced $T$ is as stated. In the associated RTP with $\text{Bern}(1/2)$ marginals for $X_{\mathbf{i}}$, the random variables $\xi_{\mathbf{i}}$ are never used, so $X_{\varnothing}$ is independent of $\mathcal{G}$ and the endogenous property fails.

*Bivariate uniqueness.*   In the setting of an RTF we have the induced map $T : \mathcal{P} \to \mathcal{P}(S)$. Now consider a bivariate version. Write $\mathcal{P}^{(2)}$ for the space of probability measures on $S^2 = S \times S$ with marginals in $\mathcal{P}$. Define $T^{(2)} : \mathcal{P}^{(2)} \to \mathcal{P}(S^2)$ by:

$T^{(2)}(\mu^{(2)})$ is the distribution of $(g(\xi, X_i^{(1)}, 1 \le i \le^* N), g(\xi, X_i^{(2)}, 1 \le i \le^* N))$, where:

  (i)  $((X_i^{(1)}, X_i^{(2)}), i \ge 1)$ are independent with distribution $\mu^{(2)}$ on $\mathcal{P}^{(2)}$;
  (ii)  $(\xi, N)$ has distribution $\nu$;
  (iii)  the families in (i) and (ii) are independent.

The point is that we use the *same realizations* of $(\xi, N)$ in both components. Immediately from the definitions we have:

  (a) If $\mu$ is a fixed point for $T$, then the associated *diagonal measure* $\mu^{\nearrow}$ is a fixed point for $T^{(2)}$, where

$$\mu^{\nearrow} = \text{dist}(X, X) \qquad \text{for } \mu = \text{dist}(X).$$

  (b) If $\mu^{(2)}$ is a fixed point for $T^{(2)}$, then each marginal distribution is a fixed point for $T$.

So if $\mu$ is a fixed point for $T$, then $\mu^{\nearrow}$ is a fixed point for $T^{(2)}$ and there may or may not be other fixed points of $T^{(2)}$ with marginals $\mu$.

Definition 9. An invariant RTP with marginal $\mu$ has the *bivariate uniqueness* property if $\mu^{\nearrow}$ is the unique fixed point of $T^{(2)}$ with marginals $\mu$.

The next example shows that even when $\mu$ is the *unique* fixed point of $T$, there may be fixed points of $T^{(2)}$ other than $\mu^{\nearrow}$.

Example 10. Take independent $I, \xi$ such that $I$ has Bern(1/2) distribution and $\xi$ has Bern($q$) distribution for some $0 < q < 1$. Consider the RDE

$$X \stackrel{d}{=} X_{I+1} + \xi \bmod 2; \qquad S = \{0, 1\}.$$

Here $T$ maps Bern($p$) to Bern($p'$) where $p' = p(1-q) + (1-p)q$, so that Bern(1/2) is the unique fixed point of $T$. But product measure Bern(1/2) $\times$ Bern(1/2) is a fixed point for $T^{(2)}$, and this differs from (Bern(1/2))$^{\nearrow}$.

2.5. *The equivalence theorem.* Here we state a version of the general result linking endogeny and bivariate uniqueness, without seeking minimal hypotheses. The result and proof are similar to standard results about Gibbs measures and Markov random fields (see Chapter 7 of [34]), but our emphasis is different, so it seems helpful to give a direct proof here, after a few remarks.

Theorem 11. *Suppose $S$ is a Polish space. Consider an invariant RTP with marginal distribution $\mu$.*

(a) *If the endogenous property holds, then the bivariate uniqueness property holds.*

(b) *Conversely, suppose the bivariate uniqueness property holds. If also $T^{(2)}$ is continuous with respect to weak convergence on the set of bivariate distributions with marginals $\mu$, then the endogenous property holds.*

(c) *Further, the endogenous property holds if and only if $T^{(2)^n}(\mu \otimes \mu) \stackrel{d}{\to} \mu^{\nearrow}$, where $\mu \otimes \mu$ is product measure.*

Here $T^{(2)^n}$ denotes the $n$th iterate of $T^{(2)}$. Note that in part (c) we do not need to assume continuity of $T^{(2)}$. Also (c) can be used nonrigorously to investigate endogeny via numerical or Monte Carlo methods, as will be described in Section 8.1. For the record we state:

Open Problem 12. Can the continuity hypothesis in (b) be removed?

Example 13 (*Noisy voter model on directed tree*). This example shows that the endogenous property may hold for some invariant measures while failing for others. Take $S = \{0, 1\}$ and let $\xi$ have Bern($\varepsilon$) distribution for small $\varepsilon > 0$. Consider the RDE

$$X \stackrel{d}{=} \xi + \mathbb{1}_{(X_1 + X_2 + X_3 \geq 2)} \bmod 2.$$

In words, a parent vertex adopts the majority opinion of its three children nodes, except with probability $\varepsilon$ adopting the opposite opinion. The Bern($p$) distribution is invariant iff $p$ satisfies

$$p = (1 - \varepsilon)q(p) + \varepsilon(1 - q(p)); \qquad q(p) = p^3 + 3p^2(1 - p).$$

There are three solutions $\{p_*(\varepsilon), \frac{1}{2}, 1 - p_*(\varepsilon)\}$ where $p_*(\varepsilon) \downarrow 0$ as $\varepsilon \downarrow 0$. As in Example 10, the invariant RTP with Bern($1/2$) marginal is not endogenous because the product measure is invariant for $T^{(2)}$. But the invariant RTP with Bern($p_*(\varepsilon)$) marginal is endogenous; one can check that $T^{(2)}$ is a strict contraction on the space of bivariate distributions with Bern($p_*(\varepsilon)$) marginals, and then appeal to the contraction lemma and to Theorem 11(c).

REMARKS. Theorem 21 and Corollary 26 provide other contexts where the endogenous property holds for the "fundamental" invariant measure but not for others. Contexts where the fundamental invariant measure is nonendogenous are less common: Proposition 48 is the most natural example.

PROOF OF THEOREM 11.   (a) Let $\nu$ be a fixed point of $T^{(2)}$ with marginals $\mu$. Consider a bivariate RTP $((X_{\mathbf{i}}^{(1)}, X_{\mathbf{i}}^{(2)}), \mathbf{i} \in \mathcal{T})$ with $\nu = \text{dist}(X_\varnothing^{(1)}, X_\varnothing^{(2)})$. Define $\mathcal{G}_n = \sigma((\xi_{\mathbf{i}}, N_{\mathbf{i}}) : \text{gen}(\mathbf{i}) \leq n)$ where for $\mathbf{i} = i_1 i_2 \cdots i_d$ we set $\text{gen}(\mathbf{i}) = d$, its generation. Observe that $\mathcal{G}_n \uparrow \mathcal{G}$.

Fix $\Lambda : S \to \mathbb{R}$ a bounded continuous function. Notice that from the construction of the bivariate RTP,

$$(X_\varnothing^{(1)}; (\xi_{\mathbf{i}}, N_{\mathbf{i}}), \text{ gen}(\mathbf{i}) \leq n) \stackrel{d}{=} (X_\varnothing^{(2)}; (\xi_{\mathbf{i}}, N_{\mathbf{i}}), \text{ gen}(\mathbf{i}) \leq n).$$

So

(11)  $$\mathbb{E}[\Lambda(X_\varnothing^{(1)})|\mathcal{G}_n] = \mathbb{E}[\Lambda(X_\varnothing^{(2)})|\mathcal{G}_n] \qquad \text{a.s.}$$

Now by martingale convergence

(12)  $$\mathbb{E}[\Lambda(X_\varnothing^{(1)})|\mathcal{G}_n] \xrightarrow{\text{a.s.}} \mathbb{E}[\Lambda(X_\varnothing^{(1)})|\mathcal{G}] \stackrel{\text{a.s.}}{=} \Lambda(X_\varnothing^{(1)}),$$

the last equality because of the *endogenous* assumption for the univariate RTP. Similarly,

$$\mathbb{E}[\Lambda(X_\varnothing^{(2)})|\mathcal{G}] \stackrel{\text{a.s.}}{=} \Lambda(X_\varnothing^{(2)}).$$

Thus by (11) we see that $\Lambda(X_\varnothing^{(1)}) = \Lambda(X_\varnothing^{(2)})$ a.s. Since this is true for every bounded continuous $\Lambda$ we deduce $X_\varnothing^{(1)} = X_\varnothing^{(2)}$ a.s., proving bivariate uniqueness.

(b) To prove the converse, again fix $\Lambda: S \to \mathbb{R}$ bounded continuous. Let $(X_{\mathbf{i}})$ be the invariant RTP with marginal $\mu$. Again by martingale convergence

$$\tag{13} \mathbb{E}[\Lambda(X_{\varnothing})|\mathcal{G}_n] \xrightarrow[\mathcal{L}_2]{\text{a.s.}} \mathbb{E}[\Lambda(X_{\varnothing})|\mathcal{G}].$$

Independently of $(X_{\mathbf{i}}, \xi_{\mathbf{i}}, N_{\mathbf{i}}, \ \mathbf{i} \in \mathbb{T})$, construct random variables $(V_{\mathbf{i}}, \mathbf{i} \in \mathbb{T})$ which are i.i.d. with distribution $\mu$. For $n \geq 1$, define $Y_{\mathbf{i}}^n := V_{\mathbf{i}}$ if $\text{gen}(\mathbf{i}) = n$, and then recursively define $Y_{\mathbf{i}}^n$ for $\text{gen}(\mathbf{i}) < n$ by (8) to get an invariant RTP $(Y_{\mathbf{i}}^n)$ of depth $n$. Observe that $X_{\varnothing} \overset{d}{=} Y_{\mathbf{i}}^n$. Further given $\mathcal{G}_n$, the variables $X_{\varnothing}$ and $Y_{\varnothing}^n$ are conditionally independent and identically distributed given $\mathcal{G}_n$. Now let

$$\tag{14} \sigma_n^2(\Lambda) := \|\mathbb{E}[\Lambda(X_{\varnothing})|\mathcal{G}_n] - \Lambda(X_{\varnothing})\|_2^2.$$

We calculate

$$\tag{15} \begin{aligned} \sigma_n^2(\Lambda) &= \mathbb{E}[(\Lambda(X_{\varnothing}) - \mathbb{E}[\Lambda(X_{\varnothing})|\mathcal{G}_n])^2] \\ &= \mathbb{E}[\text{var}(\Lambda(X_{\varnothing})|\mathcal{G}_n)] \\ &= \tfrac{1}{2}\mathbb{E}[(\Lambda(X_{\varnothing}) - \Lambda(Y_{\varnothing}^n))^2]. \end{aligned}$$

The last equality uses the conditional form of the fact that for any random variable $U$ one has $\text{var}(U) = \frac{1}{2}\mathbb{E}[(U_1 - U_2)^2]$, where $U_1, U_2$ are i.i.d. copies of $U$.

Now suppose we show that

$$\tag{16} (X_{\varnothing}, Y_{\varnothing}^n) \overset{d}{\to} (X^{\star}, Y^{\star}) \qquad \text{say, as } n \to \infty$$

for some limit. From the construction,

$$\begin{bmatrix} X_{\varnothing} \\ Y_{\varnothing}^{n+1} \end{bmatrix} \overset{d}{=} T^{(2)}\left(\begin{bmatrix} X_{\varnothing} \\ Y_{\varnothing}^n \end{bmatrix}\right),$$

and then the weak continuity assumption on $T^{(2)}$ implies

$$\begin{bmatrix} X^{\star} \\ Y^{\star} \end{bmatrix} \overset{d}{=} T^{(2)}\left(\begin{bmatrix} X^{\star} \\ Y^{\star} \end{bmatrix}\right).$$

Also by construction we have $X_{\varnothing} \overset{d}{=} Y_{\varnothing}^n \overset{d}{=} \mu$ for all $n \geq 1$, and hence $X^{\star} \overset{d}{=} Y^{\star} \overset{d}{=} \mu$. The bivariate uniqueness assumption now implies $X^{\star} = Y^{\star}$ a.s. Since $\Lambda$ is a bounded continuous function, (16) implies $\Lambda(X_{\varnothing}) - \Lambda(Y_{\varnothing}^n) \to 0$ a.s. and so using (15) we see that $\sigma_n^2(\Lambda) \longrightarrow 0$. Hence from (14) and (13) we conclude that $\Lambda(X_{\varnothing})$ is $\mathcal{G}$-measurable. This is true for every bounded continuous $\Lambda$, proving that $X_{\varnothing}$ is $\mathcal{G}$-measurable, as required.

Now all that remains is to show that a limit (16) exists. Fix $f: S \to \mathbb{R}$ and $h: S \to \mathbb{R}$, two bounded continuous functions. Again by martingale convergence

$$\mathbb{E}[f(X_{\varnothing})|\mathcal{G}_n] \xrightarrow[\mathcal{L}_1]{\text{a.s.}} \mathbb{E}[f(X_{\varnothing})|\mathcal{G}],$$

and similarly for $h$. So

$$\mathbb{E}[f(X_\varnothing)h(Y_\varnothing^n)] = \mathbb{E}[\mathbb{E}[f(X_\varnothing)h(Y_\varnothing^n)|\mathcal{G}_n]]$$
$$= \mathbb{E}[\mathbb{E}[f(X_\varnothing)|\mathcal{G}_n]\mathbb{E}[h(Y_\varnothing^n)|\mathcal{G}_n]],$$

the last equality because of conditional independence of $X_\varnothing$ and $Y_\varnothing^n$ given $\mathcal{G}_n$. Letting $n \to \infty$ and using the conditionally i.d. property gives

$$(17) \qquad \mathbb{E}[f(X_\varnothing)h(Y_\varnothing^n)] \longrightarrow \mathbb{E}[\mathbb{E}[f(X_\varnothing)|\mathcal{G}]\mathbb{E}[h(X_\varnothing)|\mathcal{G}]].$$

Moreover note that $X_\varnothing \overset{d}{=} Y_\varnothing^n \overset{d}{=} \mu$ and so the sequence of bivariate distributions $(X_\varnothing, Y_\varnothing^n)$ is tight. Tightness, together with convergence (17) for all bounded continuous $f$ and $h$, implies weak convergence of $(X_\varnothing, Y_\varnothing^n)$.

(c) First assume that $T^{(2)^n}(\mu \otimes \mu) \overset{d}{\to} \mu^\nearrow$. Then with the same construction as in part (b) we get that

$$(X_\varnothing, Y_\varnothing^n) \overset{d}{\to} (X_\varnothing, X_\varnothing).$$

Further recall that $\Lambda$ is bounded continuous; thus using (13), (14) and (15) we conclude that $\Lambda(X_\varnothing)$ is $\mathcal{G}$-measurable. This is true for any bounded continuous function $\Lambda$; thus $X_\varnothing$ is $\mathcal{G}$-measurable. So the RTP is endogenous.

Conversely, suppose that the RTP with marginal $\mu$ is endogenous. Let $\Lambda_1$ and $\Lambda_2$ be two bounded continuous functions. Note that the variables $(X_\varnothing, Y_\varnothing^n)$, as defined in part (b), have joint distribution $T^{(2)^n}(\mu \otimes \mu)$. Further, given $\mathcal{G}_n$, they are conditionally independent and have the same conditional law as of $X_\varnothing$ given $\mathcal{G}_n$. So

$$\mathbb{E}[\Lambda_1(X_\varnothing)\Lambda_2(Y_\varnothing^n)] = \mathbb{E}[\mathbb{E}[\Lambda_1(X_\varnothing)|\mathcal{G}_n]\mathbb{E}[\Lambda_2(X_\varnothing)|\mathcal{G}_n]]$$
$$\to \mathbb{E}[\mathbb{E}[\Lambda_1(X_\varnothing)|\mathcal{G}]\mathbb{E}[\Lambda_2(X_\varnothing)|\mathcal{G}]]$$
$$= \mathbb{E}[\Lambda_1(X_\varnothing)\Lambda_2(X_\varnothing)].$$

The convergence is by martingale convergence, and the last equality is by endogeny. So

$$T^{(2)^n}(\mu \otimes \mu) \overset{d}{=} (X_\varnothing, Y_\varnothing^n) \overset{d}{\to} (X_\varnothing, X_\varnothing) \overset{d}{=} \mu^\nearrow. \qquad \square$$

2.6. *Tree-structured coupling from the past.* The next lemma is clearly analogous to the *coupling from the past* (CFTP) technique for studying Markov chains [56]. That technique is part of a large circle of ideas (*graphical representations* in interacting particle systems [47]; *iterated random functions* [27]) for studying uniqueness of stationary distributions, and rates of convergence to stationarity, for Markov chains via sample path constructions.

Lemma 14. *Consider an RTF and write $\mathcal{T}$ for the associated Galton–Watson tree. Suppose there exists an a.s. finite subtree $\mathcal{T}_0 \subseteq \mathcal{T}$ containing $\varnothing$ such that the following property holds a.s. for each $\mathbf{i}$:*

> *If $\mathbf{i} \in \mathcal{T}_0$, then in the relation $X_{\mathbf{i}} = g(\xi_{\mathbf{i}}, X_{\mathbf{i}j}, 1 \leq j \leq^* N_{\mathbf{i}})$ the value of $X_{\mathbf{i}}$ is unchanged by changing the values of $X_{\mathbf{i}j}$ for which $\mathbf{i}j \notin \mathcal{T}_0$.*

*Then there is a unique invariant RTP and it is endogenous.*

*In particular, if $\mathcal{T}$ is a.s. finite (equivalently, if $\mathbb{E}[N] \leq 1$ and $P(N = 1) < 1$), then there is a unique invariant RTP and it is endogenous.*

Proof. Write $\mathrm{ht}(\mathcal{T}_0)$ for the height of $\mathcal{T}_0$. Fix $d$. Define $(X_{\mathbf{i}}^{(d)}, \mathrm{gen}(\mathbf{i}) = d)$ arbitrarily, and then use (8) recursively to define $(X_{\mathbf{i}}^{(d)}, \mathrm{gen}(\mathbf{i}) \leq d)$. The hypothesis implies that on the event $\{\mathrm{ht}(\mathcal{T}_0) < d\}$ the value of $X_{\varnothing}^{(d)}$ does not depend on the arbitrary choice of $(X_{\mathbf{i}}^{(d)}, \mathrm{gen}(\mathbf{i}) = d)$, and equals some $\mathcal{G}$-measurable random element. Letting $d \to \infty$ shows there exists some $\mathcal{G}$-measurable $X_{\varnothing}$ such that

$$P(X_{\varnothing}^{(d)} \neq X_{\varnothing}) \leq P(\mathrm{ht}(\mathcal{T}_0) \geq d) \to 0.$$

The same argument applied to a first-generation individual $j$ shows there exists $\mathcal{G}$-measurable $X_j$ such that

$$P(X_j^{(d)} \neq X_j) \leq P(\mathrm{ht}(\mathcal{T}_0) \geq d - 1) \to 0.$$

Use the i.i.d. structure of $((\xi_{\mathbf{i}}, N_{\mathbf{i}}), \mathbf{i} \in \mathbb{T})$ to show that $(X_j, j \geq 1)$ are independent and distributed as $X_{\varnothing}$. Then by the defining recursion

$$X_{\varnothing} = g(\xi, X_i, 1 \leq i \leq^* N)$$

and so $\mathrm{dist}(X_{\varnothing})$ is invariant. Moreover, in any invariant RTP it must be that $X_{\varnothing}$ is this same r.v., proving uniqueness. □

Example 8 shows that (stochastic) monotonicity of the induced map $T$ is not sufficient for endogeny. The next lemma shows that a stronger "pointwise monotonicity" condition on $g$ is sufficient.

Lemma 15. *Suppose $S = \mathbb{R}^+$. Suppose $g : \Theta^* \to \mathbb{R}^+$ is monotone for each $\theta$. That is, if $1 \leq n \leq \infty$ and $\mathbf{x} = (x_i, 1 \leq i \leq^* n)$ and $\mathbf{y} = (y_i, 1 \leq i \leq^* n)$ are such that $x_i \leq y_i \ \forall i$, then $g(\theta, \mathbf{x}) \leq g(\theta, \mathbf{y})$. Suppose that for fixed $\theta$ the map $\mathbf{x} \to g(\theta, \mathbf{x})$ is continuous w.r.t. increasing limits. Suppose that, for the induced map $T$, the sequence $(T^n(\delta_0), n \geq 0)$ is tight. Then $T^n(\delta_0) \to \mu$ weakly, where the limit $\mu$ is invariant and the associated invariant RTP is endogenous.*

PROOF. This $\mu$ is the *lower invariant measure* of Lemma 4. Let $(X_{\mathbf{i}})$ be the associated RTP. For each $d$ there is a depth-$d$ RTP $(X_{\mathbf{i}}^{(d)})$ such that $\operatorname{dist}(X_{\mathbf{i}}^{(d)}) = T^{d-\operatorname{gen}(\mathbf{i})}(\delta_0)$. Using the monotonicity hypothesis

$$0 \le X_{\varnothing}^{(1)} \le X_{\varnothing}^{(2)} \le \cdots \le X_{\varnothing} \qquad \text{a.s.}$$

Since $\operatorname{dist}(X_{\varnothing}^{(d)}) \to \mu$ we have $X_{\varnothing}^{(d)} \uparrow X_{\varnothing}$ a.s., and then since $X_{\varnothing}^{(d)}$ is $\mathcal{G}$-measurable we see that $X_{\varnothing}$ is $\mathcal{G}$-measurable. $\square$

2.7. *Markov chains.* Any Markov chain can be represented (distributionally) as an iterated random function $X_n = g(X_{n-1}, \xi_n)$ for i.i.d. $(\xi_n)$ and some $g$. So the stationary distributions (if any) are the solutions of

$$X \overset{d}{=} g(X, \xi).$$

This is the special case of RDEs for which $P(N = 1) = 1$. In general when we talk about RDEs we are envisaging the case where $P(N \ge 2) > 0$.

**3. The linear case.** The basic linear case is the case $g((\xi, X_i)) = \sum_{i=1}^{N} \xi_i X_i$ on $S = \mathbb{R}$. Note the $(\xi_i)$ may be dependent. This and the extension (20) have been studied quite extensively; our discussion focuses on analogies with the max-type cases later. Where the state space is $\mathbb{R}^+$, the key ideas are from [28] which assumed $N$ is nonrandom; the extensions to random $N$ (which is a frequent setting for our max-type examples) have been developed in [49, 50]. Here is a typical result (Corollaries 1.5 and 1.6 of [49]; the case of nonrandom $N$ is in [28]; minor nontriviality assumptions omitted).

THEOREM 16. *Suppose $\xi_i \ge 0$, with $\xi_i > 0$ iff $1 \le i \le N$, for some random $0 \le N < \infty$. Suppose $N$ and $\sum_i \xi_i$ have finite $(1 + \delta)$th moments, for some $\delta > 0$. Write $\rho(x) = \mathbb{E}[\sum_i \xi_i^x]$. Suppose there exists $0 < \alpha \le 1$ such that $\rho(\alpha) = 1$ and $\rho'(\alpha) \le 0$. Suppose either:*

(i) *$\alpha = 1$; or*

(ii) *the measure $\sum_i P(\log \xi_i \in \cdot)$ is not centered-lattice, that is to say, not supported on $s\mathbb{Z}$ for any real $s > 0$.*

*Then the RDE*

$$(18) \qquad X \overset{d}{=} \sum_i \xi_i X_i \qquad (S = \mathbb{R}^+)$$

*has an invariant distribution $X$ with $P(X = 0) < 1$, and this solution is unique up to multiplicative constants. In case* (i), *$\mathbb{E}(X) < \infty$ if $\rho'(\alpha) < 0$. In case* (ii), *if $\alpha < 1$, then $P(X > x) \sim c x^{-\alpha}$ as $x \to \infty$, for some $0 < c < \infty$.*

One can study ([49], Theorem 6.1) the operator $T$ corresponding to (18) with respect to the metric $d_\alpha$ defined as at (7) but without the $(\cdot)^{1/\alpha}$ term. In the setting of Theorem 16 it turns out that the contraction coefficient is $\rho(\alpha) = 1$ and hence the contraction argument cannot be used directly. The proof of Theorem 16 instead involves somewhat intricate analysis to find the moment generating function of $X$. See [38] for the case where $N$ may be infinite.

See [49] for many references to the appearance of the linear RDE (18) in the study of branching processes and branching random walks, invariant measures of infinite particle systems, and Hausdorff dimension of random Cantor-type sets. See [26, 60] for many references to linear RDEs arising in probabilistic analysis of algorithms which are analyzable by contraction. See [39, 40] for the specialization

$$X \overset{d}{=} \sum_{i=1}^{\infty} h(\xi_i) X_i \qquad (S = \mathbb{R}^+)$$

where $(\xi_i)$ are the points of a Poisson process on $(0, \infty)$. Often, within one model there are different questions which lead to both linear and max-type RDEs; instances can be found in Sections 4.1, 5 and 7.4.

Questions of endogeny have apparently not been studied in this linear case. Note that Example 2 provides a (degenerate?) case of a linear RDE on $\mathbb{R}$ which is not endogenous. The following corollary deals with the simplest specialization of the Theorem 16 setting.

COROLLARY 17. *In the setting of Theorem 16, suppose* (i) *holds and $\rho'(1) < 0$, so that the RDE* (18) *has a solution $X$ with $\mathbb{E}X < \infty$ and $P(X = 0) < 1$. Then the associated RTP is endogenous.*

PROOF. Consider a solution of the bivariate fixed point equation

$$(X, Y) \overset{d}{=} \left( \sum_i \xi_i X_i, \sum_i \xi_i Y_i \right).$$

Observe

$$|X - Y| \overset{d}{=} \left| \sum_i \xi_i (X_i - Y_i) \right| \leq \sum_i \xi_i |X_i - Y_i|$$

and the expectations of the leftmost and rightmost terms are equal. So the inequality must be the a.s. equality

$$\left| \sum_i \xi_i (X_i - Y_i) \right| = \sum_i \xi_i |X_i - Y_i| \qquad \text{a.s.}$$  (19)

By Theorem 11(b) it is enough to show $X = Y$ a.s. Suppose not. Then $X_i - Y_i$ takes both positive and negative values. So we cannot have $P(\xi_1 > 0, \xi_2 > 0) > 0$ or there would be nonzero chance of cancellation in the sum and (19) would fail. Thus the RDE can only be of the form $X \overset{d}{=} \xi_1 X_1$. But this can only happen if $P(\xi_1 = 1) = 1$, which case is excluded by the hypothesis $\rho'(1) < 0$.  □

OPEN PROBLEM 18. *Study endogeny in the other cases of Theorem 16.*

It is worth pointing out that there is no very complete "general theory" for $S = \mathbb{R}$:

OPEN PROBLEM 19. *Study analogs of Theorem 16 for $S = \mathbb{R}$.*

Of course the contraction method remains useful in particular cases. See [22] for results on smoothness of solutions in the case of finite second moment.

3.1. *The Quicksort RDE.* A slight extension of the linear case is the case

$$(20) \qquad g((\xi_i, X_i)) \overset{d}{=} \xi_0 + \sum_{i \geq 1} \xi_i X_i.$$

As a well-known concrete example, probabilistic analysis of the asymptotic distribution of the number of comparisons in the Quicksort algorithm leads to the study of the following RDE:

$$(21) \qquad X \overset{d}{=} U X_1 + (1 - U) X_2 + C(U) \qquad (S = \mathbb{R})$$

where $C(x) := 2x \log x + 2(1 - x) \log(1 - x) + 1$, and $U \overset{d}{=} U(0, 1)$. There is a unique solution with $\mathbb{E}[X^2] < \infty$ because $T$ is a contraction under the metric $d_2$ at (7) [59]. But there are also other solutions.

THEOREM 20 ([31]). *Let $\nu$ be the solution of the RDE (21) with zero mean and finite variance. Then the set of all solutions is the set of distributions of the form $\nu * \text{Cauchy}(m, \sigma^2)$ where $m \in \mathbb{R}$ and $\sigma^2 \geq 0$, and $*$ denotes convolution.*

The next result basically says that none other than the "fundamental" solution of the Quicksort RDE (21) is endogenous.

THEOREM 21. *An invariant RTP associated with the Quicksort RDE (21) is endogenous if and only if $\sigma = 0$.*

PROOF. Let $\mu$ be a solution of the RDE (21), so using Theorem 20 $\mu = \nu * \text{Cauchy}(m, \sigma^2)$ for some $m \in \mathbb{R}$ and $\sigma^2 \geq 0$. Suppose $(X, Y)$ is a solution of the bivariate RDE with marginals $\mu$

$$(22) \qquad \begin{pmatrix} X \\ Y \end{pmatrix} = \begin{pmatrix} UX_1 + (1 - U)X_2 + C(U) \\ UY_1 + (1 - U)Y_2 + C(U) \end{pmatrix},$$

where $(X_1, Y_1)$ and $(X_2, Y_2)$ are i.i.d. having the same distribution as $(X, Y)$ and are independent of $U \stackrel{d}{=} \text{Uniform}[0, 1]$.

First consider the case $\sigma = 0$. In this case both $X$ and $Y$ have finite second moment and hence so does $D = X - Y$. Naturally the distribution of $D$ satisfies the RDE

$$D \stackrel{d}{=} UD_1 + (1 - U)D_2 \qquad \text{(on } \mathbb{R}),$$

where $D_i = X_i - Y_i$, $i \in \{1, 2\}$. Easy calculation shows that $\mathbb{E}[D] = 0 = \mathbb{E}[D^2]$. Thus $X = Y$ a.s. which proves endogeny for the invariant RTP with marginal $\mu$, by using part (b) of Theorem 11.

Now consider the other case $\sigma > 0$. Let $Q$ be a random variable with distribution $\nu$ and let $(Z, W)$ be i.i.d. Cauchy$(m, \sigma^2)$, independent of $Q$. We claim that $(X, Y) = (Q + Z, Q + W)$ is a solution of the bivariate equation (22). In that case $X \neq Y$ a.s. and hence part (a) of Theorem 11 implies that the invariant RTP with marginal $\mu$ is not endogenous. □

So all that remains is to prove the claim, which will use the following lemma.

LEMMA 22. Let $(Z_1, Z_2)$ be i.i.d. Cauchy$(m, \sigma^2)$ for some $m \in \mathbb{R}$ and $\sigma^2 > 0$ and let $U \stackrel{d}{=} \text{Uniform}[0, 1]$ be independent of $(Z_1, Z_2)$. Then $V = UZ_1 + (1 - U)Z_2$ is a Cauchy$(m, \sigma^2)$ random variable which is independent of $U$.

PROOF. We will calculate the characteristic function of $V$ conditioned on $U$. Fix $t \in \mathbb{R}$; then

$$\mathbb{E}[e^{itV}|U] = \mathbb{E}[e^{iUtZ_1}e^{i(1-U)tZ_2}|U]$$
$$= \mathbb{E}[e^{iUtZ_1}|U] \times \mathbb{E}[e^{i(1-U)tZ_2}|U]$$
$$= \exp(imtU - \sigma U|t|) \times \exp(imt(1 - U) - \sigma(1 - U)|t|)$$
$$= \exp(imt - |t|)$$

as required.

Now let $(Q_1, Q_2)$ be two independent copies of $Q$ and let $(Z_1, Z_2, W_1, W_2)$ be i.i.d. Cauchy$(m, \sigma^2)$ which are independent of $(Q_1, Q_2)$. Define $X_i = Q_i +$ $Z_i$ and $Y_i = Q_i + W_i$ for $i \in \{1, 2\}$. Then $(X_1, Y_1)$ and $(X_2, Y_2)$ are two i.i.d. copies of $(X, Y)$. Trivially $UX_1 + (1 - U)X_2 + C(U) = Q' + Z'$ and $UY_1 + (1 - U)Y_2 + C(U) = Q' + W'$, where $Q' = UQ_1 + (1 - U)Q_2 + C(U)$, $Z' = UZ_1 + (1 - U)Z_2$ and $W' = UW_1 + (1 - U)W_2$. Notice that $Q' \overset{d}{=} Q$ and that by Lemma 22, $Z'$ and $W'$ are i.i.d. Cauchy$(m, \sigma^2)$ which are independent of $Q'$. Hence $(Q' + Z', Q' + W') \overset{d}{=} (X, Y)$. This proves the claim. □

3.2. *Moment recursions.* Another feature of the linear and the extended linear cases is that one can give a recursion for the moments of the solutions $X$, assuming moments exist. For instance, in case (20)

$$\mathbb{E}[X] = \mathbb{E}[\xi_0] + \left( \sum_{i \geq 1} \mathbb{E}[\xi_i] \right) \mathbb{E}[X],$$

$$\mathbb{E}[X^2] = \mathbb{E}[\xi_0^2] + \left( 2 \sum_{i \geq 1} \mathbb{E}[\xi_0 \xi_i] \right) \mathbb{E}[X]$$

$$+ \left( \sum_{i,j \geq 1, i \neq j} \mathbb{E}[\xi_i \xi_j] \right) (\mathbb{E}[X])^2 + \left( \sum_{i \geq 1} \mathbb{E}[\xi_i^2] \right) \mathbb{E}[X^2].$$

Unfortunately one does not have analogous general explicit information in our max-type setting.

**4. Simple examples of max-type RDEs.** The examples in this section are "simple" in a particular sense: one can construct an explicit solution (typically in terms of the stochastic process from which the RDE arises) without needing first to solve the fixed point equation analytically.

4.1. *Height of subcritical Galton–Watson trees.* A *Galton–Watson tree* is the family tree of a Galton–Watson branching process with offspring distribution $N$, say, and with one progenitor. Exclude as trivial the cases $P(N = 0) = 1$ and $P(N = 1) = 1$. In the (sub)critical case $\mathbb{E}[N] \leq 1$, it is well known by probabilistic arguments that the branching process becomes extinct a.s., so that the random variable

$H := \min\{g \,|\, \text{no individuals in generation } g\} = 1 + (\text{height of the tree})$

is a.s. finite. By conditioning on the number $N$ of offspring of the progenitor, we see that $H$ satisfies the RDE

(23) $$H \overset{d}{=} 1 + \max(H_1, H_2, \ldots, H_N), \qquad H \in \{1, 2, 3, \ldots\},$$

where the max over an empty set equals zero. Lemma 14 shows this RDE has a unique solution and is endogenous (of course this is also easy to check directly).

This RDE (23) is a natural prototype for max-type RDEs, and the following section describes one direction of generalization.

Note that the total progeny $Z$ in the Galton–Watson tree satisfies the linear RDE in Example 1. This is one of several settings where aspects of the "typical" behavior are governed by a linear RDE while aspects of the "extreme" behavior are governed by a max-type RDE.

4.2. *Positive range of one-dimensional BRW.*   Consider a discrete-generation process in which individuals are at positions on the real line $\mathbb{R}$. In generation 0 there is one individual, at position 0. In generation 1 we see that individual's offspring; there are $N$ offspring (for random $0 \leq N \leq \infty$) at positions $\infty > \xi_1 \geq \xi_2 \geq \cdots$, the joint distribution of $(N; \xi_i, i \geq 1)$ being arbitrary subject to the moment condition (24) below. Inductively, each individual in generation $n$, at position $x$ say, has $N'$ children at positions $(x + \xi'_i)$, where the families $(N'; \xi'_i)$ are i.i.d. for different parents. This process is (discrete-time, one-dimensional) *branching random walk* (BRW). The phrase "random walk" indicates the spatial homogeneity (otherwise we would have a *branching Markov chain*). Some authors use BRW for the more special case where different siblings' displacements are independent of each other and of $N$; we shall call this the IBRW (independent BRW) case. Write

$$m(\theta) = \mathbb{E}\left[\sum_i e^{\theta \xi_i}\right].$$

The moment condition we shall assume throughout is

(24)             $\exists \theta > 0$ such that $m(\theta) < \infty$.

The underlying Galton–Watson process, obtained by ignoring spatial positions, may be subcritical, supercritical or critical (depending on $\mathbb{E}[N]$, as usual). Consider

$R_n :=$ position of rightmost individual in generation $n$

with $R_n = -\infty$ if there is no such individual. Write nonextinction for the event that the process survives forever. Standard results going back to [19] show:

PROPOSITION 23.   *If the BRW is supercritical, then there exists a constant $-\infty < \gamma < \infty$ such that $n^{-1}R_n \to \gamma$ a.s. on nonextinction. And $\gamma$ is computable as the solution of*

$$\inf_{\theta > 0} (\log m(\theta) - \gamma \theta) = 0.$$

Now consider $R := \max_{n \geq 0} R_n$, the position of the rightmost particle ever. If the process becomes extinct a.s., or in the setting of Proposition 23 with $\gamma < 0$, we clearly have $0 \leq R < \infty$ a.s. Studying $R$ generalizes the study (Section 4.1) of the height of a Galton–Watson tree (take $\xi_i = 1$), as well as the study of the rightmost position of a random walk [take $N = 1$ and compare with (5)]. Applications to queueing networks are given in [44]. Conditioning on the first-generation offspring leads to the RDE below, and Lemma 15 establishes the other assertions.

LEMMA 24. *Suppose extinction is certain, or suppose supercritical with $\gamma < 0$. Then*

$$R \stackrel{d}{=} \max\left(0, \max_i(R_i + \xi_i)\right), \qquad 0 \leq R < \infty. \tag{25}$$

*If extinction is certain, then $R$ is the* unique *solution of this RDE and the RTP is endogenous. In the supercritical case, $R$ is the lower invariant measure for the RTP, and the associated invariant RTP is endogenous.*

There is an interesting *critical scaling* question here—see Open Problem 30 later. The lemma also leaves open the question of whether there may be other invariant measures in the supercritical case. A thorough treatment of the latter question was given in [21] within a slightly more general setting, including the next result showing that *nonuniqueness* is typical.

PROPOSITION 25 ([21], Theorem 1). *Under technical hypotheses (omitted here) the RDE (25) has a one-parameter family of solutions $X(\gamma)$, $0 \leq \gamma < \infty$. Here $X(0)$ is the lower invariant measure. There exists $\alpha > 0$ such that for each $\gamma > 0$ we have $P(X(\gamma) > x) \sim c_\gamma \exp(-\alpha x)$ as $x \to \infty$, for some $0 < c_\gamma < \infty$.*

Without needing to go into the proof of Proposition 25, we can observe the following.

COROLLARY 26. *Under the hypotheses of Proposition 25, for $\gamma > 0$ the invariant RTP associated with $X(\gamma)$ is not endogenous.*

PROOF. Let $(Q_\mathbf{i}, \mathbf{i} \in \mathcal{T})$ be the associated BRW; that is, $\mathcal{T}$ is the family tree of descendants of the progenitor, and $Q_\mathbf{i}$ is the position on $\mathbb{R}$ of individual $\mathbf{i}$, with $Q_\varnothing = 0$. Fix $d$ and consider the following construction. Let $(Z_\mathbf{i}^{(d)} : \text{gen}(\mathbf{i}) = d)$ be i.i.d. with some invariant distribution. For $\mathbf{i} \in \mathcal{T}$ define

$$Y_\mathbf{i}^{(d)} = Z_\mathbf{i}^{(d)}, \qquad \text{gen}(\mathbf{i}) = d,$$

and then for $\mathrm{gen}(\mathbf{i}) = d-1, d-2, \ldots, 1, 0$ define

$$Y_{\mathbf{i}}^{(d)} = \max(0;\ Q_{\mathbf{j}} - Q_{\mathbf{i}},\ \mathrm{gen}(\mathbf{j}) < d;\ Q_{\mathbf{j}} - Q_{\mathbf{i}} + Z_{\mathbf{j}}^{(d)},\ \mathrm{gen}(\mathbf{j}) = d)$$

where $\mathbf{j}$ runs over all descendants of $\mathbf{i}$. One can check that $(Y_{\mathbf{i}}^{(d)})$ defines an invariant RTP of depth $d$. Now let $A_d$ be the event that, in the definition

$$Y_{\varnothing}^{(d)} = \max(0;\ Q_{\mathbf{j}},\ \mathrm{gen}(\mathbf{j}) < d;\ Q_{\mathbf{j}} + Z_{\mathbf{j}}^{(d)},\ \mathrm{gen}(\mathbf{j}) = d)$$

the maximum is attained by some generation-$d$ descendant. We use the following straightforward lemma whose proof is given later.

LEMMA 27. *For an r.v. $Y$ and $\delta > 0$ define*

$$\mathrm{conc}(\mathrm{dist}(Y), \delta) = \max_a P(a \leq Y \leq a + \delta).$$

*Suppose $(Z_i)$ are i.i.d. with $P(Z > x) \sim ce^{-\alpha x}$ as $x \to \infty$. Then there exists $\delta > 0$, depending only on the distribution of $Z$, such that for every countable set $(x_i)$ of reals for which $Y := \max_i(x_i + Z_i) < \infty$ a.s., we have $\mathrm{conc}(\mathrm{dist}(Y), \delta) \leq 1 - \delta$.*

On the event $A_d$ the r.v. $Y_{\varnothing}^{(d)}$ is of the form in the lemma, with the role of the $(x_i)$ played by the $\mathcal{G}_d$-measurable r.v.'s $(Q_{\mathbf{j}} : \mathrm{gen}(\mathbf{j}) = d)$, where $\mathcal{G}_d$ is the $\sigma$-field generated by the first $d$ generations of the BRW. So the lemma, together with the tail estimate in Proposition 25, implies

$$\mathrm{conc}(\mathrm{dist}(Y_{\varnothing}^{(d)} | \mathcal{G}_d),\ \delta) \leq 1 - \delta \qquad \text{on } A_d.$$

This estimate remains true for an invariant RTP $(Y_{\mathbf{i}})$ of infinite depth. If the RTP is endogenous, then $Y_{\varnothing}$ is $\mathcal{G}$-measurable, and so the conditional distributions of $Y_{\varnothing}$ given $\mathcal{G}_d$ converge as $d \to \infty$ to the unit mass at $Y_{\varnothing}$; then the inequality above implies $P(A_d) \to 0$. But $P(A_d) \to 0$ implies

$$Y_{\varnothing} = \max(Q_{\mathbf{i}} : \mathbf{i} \in \mathcal{T})$$

and so the invariant distribution is just the lower invariant distribution.  □

PROOF OF LEMMA 27. Suppose if possible the conclusion of the lemma is not true. Then for every $\delta_n \downarrow 0+$ we can find a countable collection of reals $(x_i^n)_{i \geq 1}$ such that $Y_n := \max_{i \geq 1}(x_i^n + Z_i) < \infty$ a.s. and

$$P(0 \leq Y_n \leq \delta_n) \geq 1 - \delta_n. \tag{26}$$

By assumption $P(Z > x) \sim ce^{-\alpha x}$ as $x \to \infty$, so $Y_n < \infty$ a.s. implies

$$0 < \sum_{i=1}^{\infty} e^{\alpha x_i^n} < \infty. \tag{27}$$

So in particular $x_i^n \to -\infty$ as $i \to \infty$ for every $n \geq 1$. Thus without loss of generality we can assume that $(x_i^n, \ i \geq 1)$ are in decreasing order.

Let $F$ be the distribution function of $Z$, and write $\bar{F}(\cdot) = 1 - F(\cdot)$. We calculate

$$P(0 \leq Y_n \leq \delta_n)$$
$$= 1 - P(Y_n \notin [0, \delta_n])$$
$$= 1 - P(Z_i < -x_i^n \text{ for all } i \geq 1, \text{ or } Z_i > \delta_n - x_i^n \text{ for some } i \geq 1)$$
$$\leq 1 - \prod_{i=1}^{\infty} F(-\lambda - x_i^n) - \max_{i \geq 1} \bar{F}(\delta_n - x_i^n),$$

for arbitrary fixed $\lambda > 0$. So from (26) we get

$$\tag{28} \prod_{i=1}^{\infty} F(-\lambda - x_i^n) + \max_{i \geq 1} \bar{F}(\delta_n - x_i^n) \leq \delta_n.$$

But $\max_{i \geq 1} \bar{F}(\delta_n - x_i^n) = \bar{F}(\delta_n - x_1^n)$, so using (28) we get

$$\tag{29} \lim_{n \to \infty} \bar{F}(\delta_n - x_1^n) = 0 \quad \Longrightarrow \quad x_1^n \to -\infty \quad \text{as } n \to \infty.$$

Now fix $\varepsilon > 0$. By hypothesis, there exists $M > 0$ such that

$$\tag{30} (1 - \varepsilon)ce^{-\alpha x} \leq \bar{F}(x) \leq (1 + \varepsilon)ce^{-\alpha x} \qquad \text{for all } x > M - \lambda.$$

Choose $n_0 \geq 1$ such that for all $n \geq n_0$ we have $x_1^n < -M$, and hence $x_i^n < -M$ for all $i \geq 1$. Now from (28)

$$\tag{31} \delta_n \geq \bar{F}(\delta_n - x_1^n) \geq (1 - \varepsilon)ce^{-\alpha(\delta_n - x_1^n)} \quad \Longrightarrow \quad e^{\alpha x_1^n} \leq \frac{1}{c(1 - \varepsilon)} \delta_n e^{\alpha \delta_n}.$$

Further, for any $k_n$,

$$\tag{32} \begin{aligned} \prod_{i=1}^{k_n} F(-\lambda - x_i^n) &= \prod_{i=1}^{k_n} (1 - \bar{F}(-\lambda - x_i^n)) \\ &\geq (1 - \bar{F}(-\lambda - x_1^n))^{k_n} \\ &\geq (1 - (1 + \varepsilon)ce^{\alpha \lambda} e^{\alpha x_1^n})^{k_n} \\ &\geq \left(1 - \frac{1 + \varepsilon}{1 - \varepsilon} e^{\alpha \lambda} \delta_n e^{\alpha \delta_n}\right)^{k_n}, \end{aligned}$$

where the last inequality follows from (31). Now take $k_n = \frac{1}{\sqrt{\delta_n}} \uparrow \infty$ to get

$$\liminf_{n \to \infty} \prod_{i=1}^{k_n} F(-\lambda - x_i^n) \geq 1.$$

This contradicts (28). $\quad \square$

4.3. *An algorithmic aspect of BRW.* In the setting of Section 4.2—a BRW satisfying (24)—there is an algorithmic question. Suppose we are in the supercritical case [say, $P(N \geq 1) = 1$ to avoid any chance of extinction] and suppose $\gamma > 0$. So there exist individuals at large positive positions— how do we find them? Suppose we can access data only by making queries. A query

    children of progenitor?

gets an answer

    progenitor has child G at position -0.4 and child V at position
    -0.8

and a query

    children of G ?

gets an answer

G has child GF at position 0.6 and child GJ at posi- tion -1.6. There is a natural *greedy algorithm* for finding individuals with large positive positions. At each step, look at all individuals $X$ named in previous steps for which one has not already made the query `children of` $X$?; then make this query for the individual $X$ at rightmost position. See Figure 3.

This greedy algorithm was studied in [2], in the special setting of binary IBRW which we now adopt (presumably much of what we say here holds in the general BRW setting). In analyzing the performance of the greedy algorithm, a key role is played by the position $L$ of the leftmost individual ever queried. On the one hand this is given by

$$(33) \qquad L = \sup_{(w_i)} \inf_i Q_{w_i}$$

where the sup is over all lines of descent $(w_i)$ and where $Q_{\mathbf{i}}$ is the position of individual $\mathbf{i}$. On the other hand, by conditioning on the positions $(\xi_i)$ of

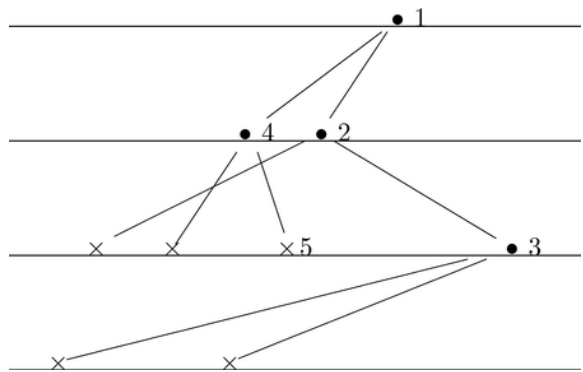

FIG. 3. *Algorithmic exploration of BRW. The individuals* ● *have been queried in order* 1, 2, 3, 4; *the children* × *have not yet been queried. Individual* 5 *will be queried next.*

the first-generation children we see that $L$ must satisfy the RDE

$$(34) \qquad L \stackrel{d}{=} \min\left(0, \max_i (L_i + \xi_i)\right), \qquad -\infty < L \le 0.$$

This is genuinely different from (25), that is, one cannot obtain (25) from (34) by, for example, reflection. As shown in [2], Proposition 4.1, this RDE has a unique solution $L$, and from (33) we see the associated invariant RTP is endogenous.

The actual question of interest in this setting in the *speed* of the greedy algorithm, defined as the limit

$$\text{speed} := \lim_n \frac{1}{n} Q_{v_n}$$

where $v_n$ is the $n$th vertex examined by the greedy algorithm. It turns out there is a simple formula for speed.

PROPOSITION 28 ([2]).  *For binary IBRW,*

$$\text{speed} = \mathbb{E}[(\xi + L)^+]$$

*where $L$ is the unique solution of the RDE (34) and $\xi$ independent of $L$ is the displacement of a child.*

So Proposition 28 is a prototype for one kind of indirect use of an RDE. From this formula, but not a priori, one can deduce that the speed is strictly positive whenever $\gamma > 0$. This leads to the question of *near-critical scaling*. Given a one-parameter family of distributions for offspring displacement $(\xi_i)$, parametrized by $p$ say, and such that $\gamma(p) > 0$ iff $p > p_{\text{crit}}$, we will typically have linear scaling for $\gamma$:

$$\gamma(p) \sim c(p - p_{\text{crit}}) \qquad \text{as } p \downarrow p_{\text{crit}}.$$

But how does $\text{speed}(p)$ scale? A special case permits explicit analysis.

THEOREM 29 ([3], Theorem 6).  *Consider an IBRW where each individual has exactly two children whose displacements $\xi$ satisfy $P(\xi = 1) = p$, $P(\xi = -1) = 1 - p$. The critical point $p_{\text{crit}}$ is the smaller solution of $16p_{\text{crit}}(1 - p_{\text{crit}}) = 1$. The solution $L(p)$ of the RDE (34) satisfies*

$$-\log P(L(p) = 0) = c(p - p_{\text{crit}})^{-1/2} + O(1) \qquad \text{as } p \downarrow p_{\text{crit}}$$

*for a certain explicitly defined constant $c$, and*

$$\text{speed}(p) = \exp(-(c + o(1))(p - p_{\text{crit}})^{-1/2}) \qquad \text{as } p \downarrow p_{\text{crit}}.$$

OPEN PROBLEM 30.  In the context of more general one-parameter families of offspring displacements $(\xi_i)$:

(a) In the supercritical setting $p \downarrow p_{\text{crit}}$, study whether the scaling for speed($p$) is as in Theorem 29.

(b) In the subcritical setting $p \uparrow p_{\text{crit}}$, study the scaling of the range $R(p)$ given in Lemma 24.

4.4. *Discounted tree sums.* In this section we study the RDE

$$(35) \qquad X \stackrel{d}{=} \eta + \max_{1 \leq i < \infty} \xi_i X_i \qquad (S = \mathbb{R}^+)$$

where $(\eta; \xi_i, 1 \leq i \leq^* N)$ has a given joint distribution, for random $N \leq \infty$. There is a natural construction of a potential solution via what we will call *discounted tree sums*, as follows. Take the associated Galton–Watson tree $\mathcal{T}$ with offspring distribution $N$. Put i.i.d. copies $\eta_{\mathbf{i}}$ of $\eta$ at vertices $\mathbf{i}$. On the edges from each $\mathbf{i}$ to its children $(\mathbf{i}j)_{j \geq 1}$ put independent copies $(\xi_{\mathbf{i}j}, 1 \leq j \leq^* N_{\mathbf{i}})$ of $(\xi_j, 1 \leq j \leq^* N)$. For an edge $e = (\mathbf{i}, \mathbf{i}j)$ we will write $\xi_e$ to denote the edge weight $\xi_{\mathbf{i}j}$. Consider a path $(\varnothing = v_0, v_1, \ldots, v_d)$. View the random variable $\eta_{v_d}$ as having "influence" $\eta_{v_d} \prod_{j=1}^{d} \xi_{(v_{j-1}, v_j)}$ at the root; that is, the influence is decreased by a factor $\xi$ in crossing an edge. From an infinite path $\pi = (\varnothing = v_0, v_1, v_2, \ldots)$ we get a total influence $\sum_{d=0}^{\infty} \eta_{v_d} \prod_{j=1}^{d} \xi_{(v_{j-1}, v_j)}$, which we suppose to be a.s. finite. Finally set

$$(36) \qquad X = \sup_{\pi = (\varnothing = v_0, v_1, v_2, \ldots)} \sum_{d=0}^{\infty} \eta_{v_d} \prod_{j=1}^{d} \xi_{(v_{j-1}, v_j)}.$$

If $X < \infty$ a.s., then clearly it is a solution of the RDE (35), and this solution is endogenous.

But it is not so easy to tell, directly from the representation (36), whether $X$ is indeed finite. So for the record we state

OPEN PROBLEM 31. *Study conditions under which* (36) *defines an a.s. finite random variable* $X$.

This question makes sense when we allow $P(\xi_i > 1) > 0$, though the concrete examples we know involve only the case $\xi_i < 1$ a.s. We content ourselves with recording a simple contraction argument (essentially that of [58], (9.1.18), in the setting of finite nonrandom $N$) designed to handle the case $\xi_i < 1$ a.s.

THEOREM 32. *Suppose* $0 \leq \xi_i < 1$ *and* $\eta \geq 0$ *with* $\mathbb{E}[\eta^p] < \infty$ $\forall p \geq 1$. *For* $1 \leq p < \infty$ *write* $c(p) := \sum_{i=0}^{\infty} \mathbb{E}[\xi_i^p] \leq \infty$. *Suppose* $c(p) < \infty$ *for some* $1 \leq p < \infty$.

(a) *The distribution* $\mu$ *of* $X$ *at* (36) *is an endogenous solution of the RDE* (35) *with all moments finite. For the associated operator* $T$ *we have* $T^n(\delta_0) \stackrel{d}{\to} \mu$.

(b) *Take $p < \infty$ such that $c(p) < 1$. Then $T$ is a strict contraction on the usual space $\mathcal{F}_p$ of distributions with finite $p$th moment. So $\mu$ is the unique solution of the RDE with finite $p$th moment, and $T^n(\mu_0) \xrightarrow{d} \mu$ for any $\mu_0 \in \mathcal{F}_p$.*

PROOF. By assumption $c(p_0) < \infty$ for some $p_0$, and then since $\xi_1 < 1$ we clearly have $c(p) \downarrow 0$ as $p \uparrow \infty$. So choose and fix $1 \leq p < \infty$ such that $c(p) < 1$. Write $\mathcal{F}_p$ for the space of distributions on $\mathbb{R}^+$ with finite $p$th moment. We will first check that $T(\mathcal{F}_p) \subseteq \mathcal{F}_p$. Let $\mu \in \mathcal{F}_p$ and let $(X)_{i \geq 1}$ be i.i.d. samples from $\mu$ which are independent of $(\xi_i)_{i \geq 1}$ and $\eta$. Define $[\mu]_p$ as the $p$th moment of $\mu$. Observe that

$$\mathbb{E}\left[\left(\max_{i \geq 1}(\xi_i X_i)\right)^p\right] = \mathbb{E}\left[\max_{i \geq 1}(\xi_i^p X_i^p)\right]$$

$$\leq \mathbb{E}\left[\sum_{i=1}^{\infty} \xi_i^p X_i^p\right]$$

$$= \sum_{i=1}^{\infty} \mathbb{E}[\xi_i^p]\mathbb{E}[X_i^p]$$

$$= [\mu]_p \times c(p) < \infty.$$

Further we have assumed that $\mathbb{E}[\eta^p] < \infty$; thus using (35) we conclude that $T$ maps $\mathcal{F}_p$ to itself.

Let $d_p$ be the usual metric (7) on $\mathcal{F}_p$. Fix $\mu, \nu \in \mathcal{F}_p$. By a standard coupling argument construct i.i.d. samples $((X_i, Y_i))_{i \geq 1}$ such that:

- they are independent of $(\xi_i)_{i \geq 1}$ and $\eta$;
- $X_i \overset{d}{=} \mu$ and $Y_i \overset{d}{=} \nu$ for all $i \geq 1$;
- $(d_p(\mu, \nu))^p = \mathbb{E}[|X_i - Y_i|^p]$.

Put $Z = \eta + \max_{i \geq 1}(\xi_i X_i)$ and $W = \eta + \max_{i \geq 1}(\xi_i Y_i)$. Notice that from definition $Z \overset{d}{=} T(\mu)$ and $W \overset{d}{=} T(\nu)$. Now

$$(d_p(T(\mu), T(\nu)))^p \leq \mathbb{E}[|Z - W|^p]$$

$$= \mathbb{E}\left[\left|\max_{i \geq 1} \xi_i X_i - \max_{i \geq 1} \xi_i Y_i\right|^p\right]$$

$$\leq \mathbb{E}\left[\sum_{i=1}^{\infty} |\xi_i X_i - \xi_i Y_i|^p\right]$$

$$= \sum_{i=1}^{\infty} \mathbb{E}[\xi_i^p] d_p^p(\mu, \nu)$$

$$= c(p) \times d_p^p(\mu, \nu).$$

So $T$ is a strict contraction map with contraction factor $(c(p))^{1/p} < 1$. Since $d_p$ is a complete metric on $\mathcal{F}_p$, the contraction method (Lemma 5) shows that there exists a fixed point $\mu \in \mathcal{F}_p$ and that $T^n(\mu_0) \xrightarrow{d} \mu$ for each $u_0 \in \mathcal{F}_p$. In particular, $T^n(\delta_0) \xrightarrow{d} \mu$. But $T^n(\delta_0)$ is just the distribution of

$$(37) \qquad X^{(n)} := \sup_{\pi=(\varnothing=v_0,v_1,v_2,\ldots,v_{n-1})} \sum_{d=0}^{n-1} \eta_{v_d} \prod_{j=1}^{d} \xi_{(v_{j-1},v_j)}.$$

So $\mu$ is the distribution of $X$ at (36). Finally, we can choose $p$ arbitrarily large, so $\mu$ has all moments finite.  $\square$

While the argument in Theorem 32, bounding a *max* by a *sum*, is crude, it does serve to establish existence of solutions in the examples we will consider below. Let us first say something about uniqueness.

COROLLARY 33.  *Consider the RDE* (35). *Suppose $X$ at* (36) *is well defined (which in particular holds under the hypotheses of Theorem 32). Write $\mu = \mathrm{dist}(X)$ for the lower invariant measure for the associated operator $T$. Consider the RDE obtained by omitting $\eta$ in* (35)*:*

$$(38) \qquad X \stackrel{d}{=} \max_{1 \le i < \infty} \xi_i X_i \qquad (S = \mathbb{R}^+).$$

*If* $\mathrm{dist}(Y)$ *is a nonzero solution of* (38)*, then for $0 \le a < \infty$ we have*

$$T^n(\mathrm{dist}(aY)) \xrightarrow{d} \mu_a$$

*and each $\mu_a$ is an invariant measure for $T$. If also*

$$(39) \qquad \eta \text{ is independent of } (\xi_i); \qquad 0 \text{ is in the support of } \eta,$$

*then each $\mu_a$ is distinct.*

By analogy with Corollary 26 and Proposition 48 later we state:

CONJECTURE 34.  *Under the assumptions of Corollary 33 and* (39)*, for $a > 0$ the invariant RTP associated with $\mu_a$ is not endogenous.*

PROOF OF COROLLARY 33.  Write $W$ for the operator associated with (38). Suppose $\nu \ne \delta_0$ is invariant for $W$. Set up the RTP $(Y_{\mathbf{i}})$ associated with (38). So for fixed $n$

$$Y_\varnothing = \sup_{\pi=(\varnothing=v_0,v_1,v_2,\ldots,v_n)} Y_{v_n} \prod_{j=1}^{n} \xi_{(v_{j-1},v_j)}$$

and the $(Y_{v_n} : \operatorname{gen}(v_n) = n)$ are independent with distribution $\nu$. And $T^n(\nu)$ is the distribution of

$$(40) \quad Z^{(n)} = \sup_{\pi = (\varnothing = v_0, v_1, v_2, \ldots, v_n)} \left( Y_{v_n} \prod_{j=1}^{n} \xi_{(v_{j-1}, v_j)} + \sum_{d=0}^{n-1} \eta_{v_d} \prod_{j=1}^{d} \xi_{(v_{j-1}, v_j)} \right).$$

Note the sample path monotonicity property

$$Z^{(n)} \le Z^{(n+1)} \qquad \text{a.s.}$$

which holds for the following reason. Given $v_n$, there is a $v_{n+1} = v_n i$ attaining the maximum $Y_{v_n} = \max_i \xi_{v_n i} X_{v_n i}$, and the right-hand side of (40) for $(v_0, v_1, \ldots, v_{n+1})$ is not smaller than the right-hand side of (40) for $(v_0, v_1, \ldots, v_n)$.

This monotonicity, together with the facts

$$Z^{(n)} \le X^{(n)} + Y_\varnothing; \qquad X^{(n)} \uparrow X \text{ a.s.},$$

implies existence of the limit

$$Z^{(n)} \uparrow Z < \infty \qquad \text{a.s.}$$

So

$$T^n(\nu) = \operatorname{dist}(Z^{(n)}) \to \operatorname{dist}(Z) := \mu_1 \qquad \text{say}$$

and by continuity, $\mu_1$ is invariant for $T$.

Now for arbitrary $a \ge 0$ define $Z_a^{(n)}$ by replacing $Y_{v_n}$ by $a Y_{v_n}$ in (40). As above there exists a limit $Z_a^{(n)} \uparrow Z_a < \infty$ a.s., and $\mu_a := \operatorname{dist}(Z_a)$ is invariant for $T$.

To prove the final assertion of Corollary 33, fix $0 < a < b$. Clearly $Z_a \le Z_b$ a.s., so it is enough to prove $P(Z_b > Z_a) > 0$. By Lemma 35 below (whose easy proof we omit) it is enough to prove

$$(41) \qquad P(a Y_\varnothing > X) > 0.$$

Let $\mathcal{H}$ be the $\sigma$-field generated by the RTP $(Y_{\mathbf{i}})$ and by all the $\xi_{v,v'}$. By assumption (39) the $(\eta_{\mathbf{i}})$ are independent of $\mathcal{H}$, and it easily follows that

$$P(X \le \varepsilon | \mathcal{H}) > 0 \qquad \text{a.s., for each } \varepsilon > 0.$$

Since $Y_\varnothing$ is $\mathcal{H}$-measurable, this establishes (41). $\square$

LEMMA 35. *For $i = 1, 2$ let $f_i \ge 0$ be a function such that $f_i^* := \sup f_i < \infty$. For $a \ge 0$ let $q(a) := \sup(a f_1 + f_2)$. If $a f_1^* > f_2^*$, then $q(b) > q(a)$ for all $b > a$.*

Corollary 33 hints that general solutions of the "discounted tree sum" RDE (35) correspond to solutions of the homogeneous RDE (38). Unfortunately the latter is not trivial to solve. For by taking logs (set $\hat{X} = \log X$, $\hat{\xi} = \log \xi$) we see (38) is equivalent to

$$\hat{X} \stackrel{d}{=} \max_i (\hat{\xi}_i + \hat{X}_i)$$

and this RDE, to be studied in Section 5, is the fundamental example of a max-type RDE which cannot be solved by any simple probabilistic construction. See [41] for further discussion of (38).

4.5. *Examples of discounted tree sums.*

EXAMPLE 36.   Take $U \stackrel{d}{=} \text{Uniform}(0,1)$ and consider the RDE

$$(42) \qquad X \stackrel{d}{=} 1 + \max(UX_1, (1-U)X_2) \qquad (S = \mathbb{R}^+).$$

This arises [25] in the context of the probabilistic worst-case analysis of Hoare's FIND algorithm. Theorem 32 implies existence of a fixed point with all moments finite, unique amongst possible fixed points with finite $(1+\varepsilon)$th moments. In a different way it can be proved [25] that any fixed point has all moments finite, and hence the fixed point is unique.

EXAMPLE 37.   Consider the RDE

$$(43) \qquad X \stackrel{d}{=} \eta + c \max(X_1, X_2)$$

where $0 < c < 1$.

This arises [13] as a "discounted branching random walk." One interpretation is as nonhomogeneous percolation on the planted binary tree (the root has degree 1), where an edge at depth $d$ has traversal time distributed at $c^d \eta$. Then $X$ is the time for the entire tree to be traversed. Assuming $\eta$ has all moments finite, Theorem 32 implies existence of a fixed point with all moments finite, unique amongst possible fixed points with finite expectation. The same conclusion can be drawn in the slightly more general setup of a Galton–Watson branching tree with offspring distribution $N$. Instead of assuming $\eta$ has all moments finite, make the weaker assumption that there exists $\theta > 0$ such that $\sup_{x \in \mathbb{R}} x^\theta P(\eta > x) < \infty$ and $mc^\theta < 1$ where $m = \mathbb{E}[N]$. Under these assumptions, [13] shows that the RDE (43) has a solution such that $P(X > x) = o(x^{-\alpha})$ where $\alpha = -\log m / \log c$. Moreover, this solution is unique in the class of distributions $H$ such that $x^\alpha(1 - H(x)) \to 0$ as $x \to \infty$. But outside this class there may be other solutions.

EXAMPLE 38. Consider the RDE

$$(44) \qquad X \stackrel{d}{=} \eta + \max_{i \geq 1} e^{-\xi_i} X_i \qquad (S = \mathbb{R}^+),$$

where $(\xi_i, i \geq 1)$ are the points of a Poisson rate 1 process on $(0, \infty)$ and where $\eta$ has Exponential(1) distribution independent of $(\xi_i)$.

This is a new example, arising from a species competition model [29]. Time reversal of the process in [29], together with a transformation of $(0, 1)$ to $(0, \infty)$ by $x \to -\log(1 - x)$, yields a branching Markov process taking values in the space of countable subsets of $(0, \infty)$, which can then be extended to $(-\infty, \infty)$ as follows.

Each individual at position $x$ at time $t$ lives for an independent Exponential($e^x$) lifetime, after which it dies and instantaneously gives birth to an infinite number of children to be placed at positions $(x + \xi_i)_{i \geq 1}$ where $(\xi_i)_{i \geq 1}$ are points of an independent Poisson rate 1 process on $(0, \infty)$. The result of [29], transformed as above, shows that for each $\lambda < \infty$ the Poisson rate $\lambda$ process on $(-\infty, \infty)$ is a stationary law for this branching Markov process. We pose a different question. What is the extinction time $X$ for the process started at time 0 with a single particle at position 0? It is easy to see that $X$ satisfies the RDE (44).

For this example, easy calculation shows that $c(p) = 1/p$ for $p \geq 1$, so Theorem 32 implies existence of an invariant distribution with all moments finite which is also endogenous.

Now in the setting of Corollary 33 consider the *homogeneous* equation, that is, with $\eta \equiv 0$

$$(45) \qquad X \stackrel{d}{=} \max_{i \geq 1} e^{-\xi_i} X_i \qquad \text{on } S = \mathbb{R}^+.$$

The solution $X \equiv 0$ of (45) corresponds to the solution of (44) with all moments finite. We show below by direct calculation that there are other solutions of (45), which by Corollary 33 correspond to other solutions of (44) with infinite mean.

PROPOSITION 39. *The set $(X_a, \ a \geq 0)$ of all solutions of the RDE* (45) *is given by*

$$(46) \qquad P(X_a \leq x) = \begin{cases} 0, & \text{if } x < 0, \\ \dfrac{x}{a + x}, & \text{if } x \geq 0. \end{cases}$$

*In particular for $a = 0$ it is the solution $\delta_0$.*

Proof. Let $\mu$ be a solution of (45). Notice that the points $\{(\xi_i; X_i)|i \geq 1\}$ form a Poisson point process, say $\mathfrak{P}$, on $(0, \infty)^2$ with mean intensity $dt\,\mu(dx)$. Thus if $F(x) = P(X \leq x)$, then for $x > 0$

$$
\begin{aligned}
F(x) &= P\left(\max_{i \geq 1} e^{-\xi_i} X_i \leq x\right)\\
&= P(\text{no points of } \mathfrak{P} \text{ are in } \{(t, z)|e^{-t}z > x\})\\
&= \exp\left(-\iint_{e^{-t}z > x} dt\,\mu(dx)\right)\\
&= \exp\left(-\int_x^\infty \frac{1 - F(u)}{u}\,du\right).
\end{aligned}
$$

(47)

We note that $F$ is infinitely differentiable so by differentiating (47) we get

(48) $$\frac{dF}{dx} = \frac{F(x)(1 - F(x))}{x} \qquad \text{for } x > 0.$$

It is easy to solve (48) to verify that the set of all solutions is given by (46). □

Later results [Proposition 48(a) applied after taking logarithms] imply that for $a > 0$ the invariant RTP associated with $X_a$ is not endogenous.

EXAMPLE 40. Nonhomogeneous percolation on the binary tree.

One can also consider the analog of (35) when $max$ is replaced by $min$, though this situation does not seem to have been studied generally. One particular occurrence is in the setting of Example 37, interpreted as non-homogeneous percolation, in which case the time $X$ taken to percolate to infinity satisfies the RDE

(49) $$X \stackrel{d}{=} \eta + c\min(X_1, X_2).$$

This setting has been studied from a different viewpoint in [17].

4.6. *Matchings on Galton–Watson trees.* Amongst many possible examples involving Galton–Watson trees, the following rather subtle example provides a warm-up to the harder example in Section 7.3.

Consider an a.s. finite Galton–Watson tree $\mathcal{T}$ with offspring distribution $N$. Fix an arbitrary probability distribution $\nu$ on $(0, \infty)$. Attach independent $\nu$-distributed weights to the edges. A *partial matching* on $\mathcal{T}$ is a subset of edges such that no vertex is in more than one edge. The *weight* of a partial matching is the sum of its edge weights. So associated with the random tree $\mathcal{T}$ is a random variable

$$W := \text{maximum weight of a partial matching.}$$

In seeking to study $W$ via recursive methods, we quickly realize that a more tractable quantity to study is

(50) $\quad X :=$ maximum weight of a partial matching
$\quad\quad\quad - $ maximum weight of a partial matching
$\quad\quad\quad\quad$ which does not include the root.

To see why, fix a child $i$ of the root. Compare (a) the maximum-weight partial matching $\mathcal{M}_i$ which matches the root to $i$ with (b) the maximum-weight partial matching $\mathcal{M}_-$ in which the root is not matched.

These matchings must agree on the subtrees of all first-generation children except $i$. On the subtree rooted at $i$, $\mathcal{M}_i$ is the maximum-weight partial matching which does not include $i$, and $\mathcal{M}_-$ is the maximum-weight partial matching. Thus weight($\mathcal{M}_i$) − weight($\mathcal{M}_-$) = $\xi_i - X_i$ where $\xi_i$ is the weight on edge (root, $i$) and $X_i$ is defined as $X$ but in terms of the subtree rooted at $i$. Since in seeking the maximum-weight partial matching we can use any $i$, or no $i$, we deduce the RDE

(51) $$X \stackrel{d}{=} \max(0, \xi_i - X_i, 1 \le i \le N)$$

where the $X_i$ are independent copies of $X$ and the $\xi_i$ are i.i.d. $(\nu)$. Uniqueness and endogeny follow from Lemma 14.

This RDE, in the special case where $N$ has Poisson(1) distribution, arose in the context of the problem

study the maximum weight $W_n$ of a partial matching on a uniform
random $n$-vertex tree, in the $n \to \infty$ limit.

The essential idea is that a randomly chosen edge of that tree splits it into two subtrees, the smaller of which is distributed as a Galton–Watson tree with Poisson(1) offspring. For the detailed story see Section 3 of [11], whose final result is:

THEOREM 41. *Suppose $\nu$ is nonatomic and has finite mean. Then*

$$\lim_n n^{-1} \mathbb{E} W_n = \mathbb{E} \xi \mathbb{1}_{(\xi > X + Z)}$$

*where the r.v.'s on the right are independent, $\xi$ has distribution $\nu$, $X$ is distributed as the solution of the RDE* (51) *with* Poisson(1) *distributed $N$ and $Z$ is distributed as the solution of the RDE*

$$Z \stackrel{d}{=} \max(X, \xi - Z)$$

*where the r.v.'s on the right are independent.*

Let us mention the explicit solutions of our RDE in two special cases.

LEMMA 42. *Let $N$ have* Poisson(1) *distribution.*

(a) ([11], *Lemma 2) If $\nu$ is the* exponential(1) *distribution, then the solution of the RDE* (51) *is*

$$P(X \le x) = \exp(-ce^{-x}), \qquad x \ge 0,$$

*where $c \approx 0.715$ is the unique strictly positive solution of $c^2 + e^{-c} = 1$.*

(b) *If $\nu$ is the* Bern$(p)$ *distribution, then the solution of the RDE* (51) *is the* Bern$(1 - x(p))$ *distribution, where $x = x(p)$ solves $x = e^{-px}$.*

A closely related "dual" problem concerns independent sets. Recall that an *independent set* in a graph is a subset of vertices, no two of which are linked by an edge. Take as before a Galton–Watson tree with $N$ offspring, and a probability distribution $\nu$ on $(0, \infty)$. Now assign independent $\nu$-distributed random weights to each *vertex* and consider

$X :=$ maximum weight of an independent set
      $-$ maximum weight of an independentset which does not include the root.

Similar to above, we can argue that $X$ is the solution of the RDE

$$(52) \qquad X \stackrel{d}{=} \max\left(0, \xi - \sum_{i=1}^{N} X_i\right)$$

where $\xi$ has distribution $\nu$, independent of $N$.

**5. Rightmost position of BRW.** Here we work in the setting of Section 4.2. We have a BRW on $\mathbb{R}$, where an individual has a random number $N$ of offspring, whose random displacements from the parent's position are $\infty > \xi_1 \ge \xi_2 \ge \cdots$, distributed arbitrarily subject to the moment condition (24). For simplicity suppose $N \ge 1$ a.s. and $P(N > 1) > 0$. By Proposition 23, the position $R_n$ of the rightmost individual in generation $n$ satisfies

$$n^{-1} R_n \to \gamma \qquad \text{a.s.}$$

For reasons to be explained in the next section, one expects that under minor extra assumptions (including a nonlattice assumption) a much stronger result is true: there exist constants $\gamma_n$ such that

$$(53) \qquad R_n - \gamma_n \stackrel{d}{\to} X \qquad \text{as } n \to \infty$$

and that $X$ is the unique (up to translation) solution of the RDE

$$(54) \qquad X \stackrel{d}{=} -\gamma + \max_i(\xi_i + X_i), \qquad -\infty < X < \infty.$$

This is our first example of an RDE which is "not simple," in the sense that we do not know how to construct probabilistically a solution.

5.1. *Tightness of $R_n$.* At first sight it may be surprising that a limit (53) could hold, since it presupposes that the sequence $(R_n - \mathrm{median}(R_n))$ is tight, whereas one might expect its spread to increase to infinity. However, tightness is quite easy to understand.

LEMMA 43. *If*

$$(55) \qquad (\mathrm{median}(R_{n+1}) - \mathrm{median}(R_n), n \geq 0) \qquad \textit{is bounded above,}$$

*then*

$$(56) \qquad (R_n - \mathrm{median}(R_n), n \geq 1) \qquad \textit{is tight.}$$

Harry Kesten (personal communication) attributes this type of argument to old work of Hammersley: it is perhaps implicit in [36], page 662.

PROOF OF LEMMA 43. Given $\varepsilon > 0$ we can choose $k < \infty$ and $B > -\infty$ such that

$$P(\text{generation } k \text{ has at least } \log_2 1/\varepsilon \text{ individuals in } [B, \infty)) \geq 1 - \varepsilon.$$

Then by conditioning on the positions of generation $k$,

$$P(R_{n+k} < B + \mathrm{median}(R_n)) \leq 2\varepsilon.$$

Writing $A$ for an upper bound in (55), we deduce

$$P(R_{n+k} < B - Ak + \mathrm{median}(R_{n+k})) \leq 2\varepsilon.$$

This establishes the tightness requirement for the left tail of $R_n$. For the right tail, given $\varepsilon < 1/4$ we can choose $k < \infty$ and $B > -\infty$ such that

$$P(\text{generation } k \text{ has at least } \varepsilon^{-1} \log 1/\varepsilon \text{ individuals in } [B, \infty)) \geq 1 - \varepsilon$$

(we changed the $\log_2 1/\varepsilon$ term above to $\varepsilon^{-1} \log 1/\varepsilon$). Write $q_n$ for the $1 - \varepsilon$ quantile of $R_n$, so that $P(R_n \geq q_n) \geq \varepsilon$. Again by conditioning on the positions of the generation $k$,

$$P(R_{n+k} < B + q_n) \leq \varepsilon + (1 - \varepsilon)^{\varepsilon^{-1} \log 1/\varepsilon} \leq 2\varepsilon < 1/2.$$

So $\mathrm{median}(R_{n+k}) \geq B + q_n$, implying

$$q_n \leq \mathrm{median}(R_{n+k}) - B \leq \mathrm{median}(R_n) + Ak - B.$$

Since $q_n$ is the $1 - \varepsilon$ quantile of $R_n$, this establishes the tightness requirement for the right tail of $R_n$. □

In the case where all displacements $\xi_i$ are nonpositive (i.e., by reflection the case where displacements are nonnegative and we are studying the position of the *leftmost* particle) it is clear that $R_n$ is decreasing and so (55) holds automatically, and then the lemma implies (56). A slicker argument for this case is in [24], Proposition 2. The same holds (by translation) if there is a constant upper bound on displacement, that is, if

$$P(\xi_1 \le x_0) = 1 \qquad \text{for some constant } x_0 < \infty.$$

From these tightness results it is natural to expect that, under rather weak regularity conditions, we in fact have the convergence results

$$\text{(57)} \qquad \text{median}(R_{n+1}) - \text{median}(R_n) \to \gamma,$$

$$\text{(58)} \qquad R_n - \text{median}(R_n) \xrightarrow{d} X,$$

for some limit distribution $X$. Our interest in these limits arises, of course, because if (57) and (58) hold, then the limit $X$ must satisfy the RDE (54).

5.2. *Limit theorems.* This topic has been studied carefully only in the IBRW setting. We quote a recent result from [14], which provides an extensive bibliography of earlier work. The proof uses a mixture of analytic and probabilistic tools, for example, the "stretching" partial order (which goes back to the original KPP paper [46]), and multiplicative martingales.

THEOREM 44 ([14]).   *Consider an IBRW where $\mathbb{E}[N] < \infty$, $N \ge 1$, $P(N > 1) > 0$, and where the offspring displacement has density $f(x) = e^{-\kappa(x)}$ for some convex function $\kappa$. Then the limit*

$$R_n - \text{median}(R_n) \xrightarrow{d} X$$

*exists [and hence satisfies the RDE (54)]. If $\mathbb{E}[N \log N] < \infty$ and if a technical assumption on $\phi(\theta) := \mathbb{E}[N] \int e^{\theta x} f(x) \, dx$ (details omitted) holds, then the limit distribution is of the form*

$$\text{(59)} \qquad P(X \le x) = \mathbb{E}[\exp(-\exp(\theta_0(Y + x)))]$$

*for some constant $\theta_0$ and random variable $Y$.*

While the log-concave assumption plays a key role in the proof, it does not seem intuitively to be essential for the result.

OPEN PROBLEM 45.   Under what weaker hypotheses does Theorem 44 remain true?

5.3. *Endogeny.* The viewpoint of this survey is to seek to study existence and uniqueness of solutions of RDEs separately from weak convergence questions. This has not been done very systematically in the present context:

OPEN PROBLEM 46.    *Study existence and uniqueness of solutions to* (54) *directly from its definition.*

However, Proposition 48 will show that the associated RDE is generally not endogenous.

We first need to exclude a degenerate case. Write

$$\gamma^* = \operatorname{ess\,sup} \xi_1.$$

If $\gamma^* < \infty$ and $\mathbb{E}\#\{i|\xi_i = \gamma^*\} > 1$, then there exist embedded infinite Galton–Watson trees on which the parent–child displacement equals $\gamma^*$; it easily follows that there is the a.s. limit

$$R_n - n\gamma^* \to X \qquad \text{a.s.}$$

and that the associated invariant RTP is endogenous. The next lemma (whose easy proof is omitted) excludes this case.

LEMMA 47.    *Consider a BRW satisfying* (24) *and* $P(N \geq 1) = 1$. *If*

$$(60) \qquad \gamma^* = \infty; \quad or \quad \gamma^* < \infty \quad and \quad \mathbb{E}\#\{i|\xi_i = \gamma^*\} < 1,$$

*then* $n^{-1}R_n \to \gamma < \gamma^*$.

PROPOSITION 48.    *Suppose $X$ is a solution of the RDE* (54). *Under either assumption* (a) *or assumption* (b) *below, the invariant RTP associated with this solution is not endogenous.*

(a) *There exist constants $c, \alpha > 0$ such that $P(X > x) \sim c\exp(-\alpha x)$ as $x \to \infty$.*

(b) *Suppose there is a BRW satisfying* (24), (60) *with $P(N \geq 1) = 1$ and $P(N > 1) > 0$. Suppose there exist constants $\gamma_n$ such that*

$$R_n - \gamma_n \xrightarrow{d} X, \qquad \gamma_n - \gamma_{n-1} \to \gamma,$$

*so that necessarily $\gamma = \lim_n n^{-1}R_n$ and $X$ satisfies the RDE* (54).

OPEN PROBLEM 49.    *Weaken the assumptions in Proposition* 48. *In particular, does nonendogeny hold under only the assumptions of Lemma* 47?

From the viewpoint of the underlying BRW, nonendogeny is a property of the RTP associated with an $n \to \infty$ limit, so it is not obvious what its significance for the BRW might be. Informally, the argument in Section 5.4 shows that nonendogeny is related to a kind of "nonpredictability" property of $R_n$. Given the ordered positions $(X_{n,i})$ of the $n$th-generation individuals, for $N > n$ write $(p_{n,N}(X_{n,i}),\ i \geq 1)$ for the probability that the rightmost individual in generation $N$ is a descendant of the $X_{n,i}$ individual. Then there exist limits $p_{n,\infty}(X_{n,i}) = \lim_{N \to \infty} p_{n,N}(X_{n,i})$. For fixed $n$ this is maximized at the rightmost individual $X_{n,1}$, but it can be shown under suitable conditions that $p_{n,\infty}(X_{n,1}) \to 0$ as $n \to \infty$. Loosely, it is unpredictable which line of descent leads to the rightmost individual at large times.

5.4. *Proof of Proposition* 48. Using the notation of Corollary 26 let $(Q_{\mathbf{i}}, \mathbf{i} \in \mathcal{T})$ be the associated BRW; that is, $\mathcal{T}$ is the family tree of the progenitor, and $Q_{\mathbf{i}}$ is the position on $\mathbb{R}$ of the $\mathbf{i}$th individual, with $Q_{\varnothing} = 0$. Fix $d \geq 1$ and let $\{Z_{\mathbf{i}}^{(d)} | \mathrm{gen}(\mathbf{i}) = d\}$ be i.i.d. copies of $X$. For $\mathbf{i} \in \mathcal{T}$ define

- $Y_{\mathbf{i}}^{(d)} = Z_{\mathbf{i}}^{(d)}$, when $\mathrm{gen}(\mathbf{i}) = d$;
- $Y_{\mathbf{i}}^{(d)} = \max\{Q_{\mathbf{j}} - Q_{\mathbf{i}} + Z_{\mathbf{j}}^{(d)} | \mathrm{gen}(\mathbf{j}) = d$ and $\mathbf{j}$ is a descendant of $\mathbf{i}\}$, when $\mathrm{gen}(\mathbf{i}) \in \{d-1, d-2, \ldots, 1, 0\}$.

It is easy to check that $(Y_{\mathbf{i}}^{(d)})$ defines an invariant RTP of depth $d$, for the RDE (54).

Let $\mathcal{G}_d$ be the $\sigma$-field generated by the first $d$ generations of the BRW. So $\mathcal{G}_d \uparrow \mathcal{G}$, the $\sigma$-field generated by all the $\xi_{\mathbf{i}}$'s. Observe that

$$(61) \qquad Y_{\varnothing}^{(d)} = \max\{Q_{\mathbf{j}} + Z_{\mathbf{j}}^{(d)} | \mathrm{gen}(\mathbf{j}) = d\}.$$

*Case* (a). Under the conditional distribution given $\mathcal{G}_d$, the random variable $Y_{\varnothing}^{(d)}$ has the same form as in Lemma 27, with the role of the $(x_i)$ being played by the $\mathcal{G}_d$-measurable random variables $(Q_{\mathbf{j}}, \mathrm{gen}(\mathbf{j}) = d)$, and the role of the $(Z_i)$ being played by the i.i.d. random variables $(Z_{\mathbf{j}}^{(d)}, \mathrm{gen}(\mathbf{j}) = d)$. So Lemma 27 along with our assumption (a) of exponential right tail, implies that there exists $\delta > 0$ such that

$$(62) \qquad \mathrm{conc}(\mathrm{dist}(Y_{\varnothing}^{(d)} | \mathcal{G}_d), \delta) \leq 1 - \delta.$$

This inequality is true for any invariant RTP of depth at least $d$, so in particular true for the invariant RTP $(Y_{\mathbf{i}})$ of infinite depth, so we get

$$\mathrm{conc}(\mathrm{dist}(Y_{\varnothing} | \mathcal{G}_d), \delta) \leq 1 - \delta \quad \Longrightarrow \quad \max_{-\infty < a < \infty} P(a \leq Y_{\varnothing} \leq a + \delta | \mathcal{G}_d) \leq 1 - \delta.$$

Now suppose that the invariant RTP were endogenous, that is, $Y_{\varnothing}$ is $\mathcal{G}$-measurable. Using the martingale convergence theorem we get for each rational $a$

$$\mathbb{1}_{(a \leq Y_{\varnothing} \leq a + \delta)} \leq 1 - \delta \qquad \text{a.s.}$$

which is clearly impossible.

For case (b) we need two lemmas. The first is straightforward (proof omitted) and the second is an analog of Lemma 27.

LEMMA 50.  *Let $p_0 < 1$. For each $n$ let $(C_{n,i}, i \geq 1)$ be independent events with $P(C_{n,i}) \leq p_0$. Suppose there is a random variable $M^*$ taking values in $\bar{\mathbb{Z}}^+ = \{0, 1, 2, \ldots; \infty\}$ such that*

$$\sum_i \mathbb{1}_{C_{n,i}} \xrightarrow{d} M^* \qquad as \; n \to \infty$$

*in the sense of convergence in distribution on $\bar{\mathbb{Z}}^+$. Then either $P(M^* = 0) > 0$ or $P(M^* = \infty) = 1$.*

LEMMA 51.  *Let $(Z_i)$ be i.i.d. nonconstant. For each $n$ let $(a_{n,i}, i \geq 1)$ be real constants. For $k \geq 1$ let $M_{n,k}$ be the $k$th largest of $(a_{n,i} + Z_i, \; i \geq 1)$. If $M_{n,1} \xrightarrow{p} 0$ as $n \to \infty$, then for each $k$ we have $M_{n,k} \xrightarrow{p} 0$ as $n \to \infty$.*

PROOF.  Write $\theta^* = \text{ess sup } Z_i$. Arrange $(a_{n,i}, i \geq 1)$ in decreasing order. Since $a_{n,1} + Z_1$ is asymptotically not greater than 0 it is easy to see that $\limsup_n a_{n,1} \leq -\theta^*$. From nonconstancy of $Z_1$ it follows that for all $\varepsilon > 0$ there exist $p_0 < 1$ and $n_0 < \infty$ such that

$$P(a_{n,1} + Z_1 \geq -\varepsilon) \leq p_0, \qquad n \geq n_0.$$

Apply Lemma 50 to the events $\{a_{n,i} + Z_i \geq -\varepsilon\}$, passing to a subsequence to assume existence of a limit

$$\sum_i \mathbb{1}_{(a_{n,i} + Z_i \geq -\varepsilon)} \xrightarrow{d} M^*.$$

By assumption $P(M^* = 0) = 0$, so by Lemma 50 $P(M^* = \infty) = 1$, implying $M_{n,k} \xrightarrow{p} 0$.

*Case* (b). Recall the argument leading to (61). Take $(\widetilde{Z}_{\mathbf{j}}^{(d)})$ to be further i.i.d. copies of $X$ and set

$$\widetilde{Y}_{\varnothing}^{(d)} = \max\{Q_{\mathbf{j}} + \widetilde{Z}_{\mathbf{j}}^{(d)} | \text{gen}(\mathbf{j}) = d\}.$$

Then the joint distribution $(Y_{\varnothing}^{(d)}, \widetilde{Y}_{\varnothing}^{(d)})$ is the distribution $T^{(2)^n}(\mu \otimes \mu)$ appearing in Theorem 11(c), and that theorem asserted that endogeny is equivalent to

(63)  $$(Y_{\varnothing}^{(d)}, \widetilde{Y}_{\varnothing}^{(d)}) \xrightarrow{d} (X, X) \qquad as \; d \to \infty.$$

Suppose, to obtain a contradiction, that (63) were true. Writing $A_d$ for the $\mathcal{G}_d$-measurable r.v. defined as the median of the conditional distribution of $Y_{\varnothing}^{(d)}$ given $\mathcal{G}_d$, (63) would easily imply

$$Y_{\varnothing}^{(d)} - A_d \xrightarrow{p} 0.$$

Now for $k \geq 1$ consider

$$B_{d,k} = k\text{th largest of } \{Q_{\mathbf{j}} + Z_{\mathbf{j}}^{(d)} | \text{gen}(\mathbf{j}) = d\}.$$

So $B_{d,1} - A_d \xrightarrow{p} 0$. Now apply Lemma 51 conditionally on $\mathcal{G}_d$, with the role of the $(a_{n,i})$ being played by $(Q_{\mathbf{j}} - A_d)$, to conclude that for each $k$ we have $B_{d,k} - A_d \xrightarrow{p} 0$. (More pedantically, we need to detour through a subsequence argument to justify conditional application of Lemma 51; we omit details.) So

$$(64) \qquad\qquad B_{d,1} - B_{d,k} \xrightarrow{p} 0.$$

Next we exploit the underlying BRW. Write

$$R_{m,k} = \text{ position of } k\text{th rightmost individual in generation } m.$$

Fix $d$ and $u > 0$. For an individual $\mathbf{j}$ in generation $d$, the displacement of its rightmost descendant in generation $m$ is asymptotically $(m \to \infty)$ distributed as $X$, independently as $\mathbf{j}$ varies, and so

$$\liminf_m P(R_{m,1} - R_{m,k} \leq u) \geq P(B_{d,1} - B_{d,k} < u)$$

by considering the rightmost descendant of each of the $k$ generation-$d$ individuals featuring in the definition of $B_{d,k}$. Now (64) implies

$$R_{m,1} - R_{m,k} \xrightarrow{p} 0 \qquad \text{as } m \to \infty.$$

But this property (for each $k$) states that an ever-increasing number of individuals accumulate near the rightmost individual, easily implying

$$R_{m+1,1} - R_{m,1} - \text{ess sup } \xi_1 \xrightarrow{p} 0.$$

This in turn implies $\lim_m m^{-1} R_{m,1} = \text{ess sup } \xi_1$, contradicting Lemma 47. $\square$

REMARK 1.   Some multiplicative martingales used in the study of BRW (see, e.g., [20]) are of the form

$$Z_n(\theta) = m^{-n}(\theta) \sum_i \exp(\theta Y_i^n)$$

where $(Y_i^n, i \geq 1)$ are the positions of the generation-$n$ individuals. The a.s. limit $Z(\theta) = \lim_n Z_n(\theta)$ satisfies the RDE

$$Z \overset{d}{=} \sum_i \exp(\theta \xi_i) Z_i / m(\theta)$$

This is an instance of an "average-case" RDE paralleling the "extreme-case" RDE (54).

REMARK 2. Very recently Iksanov has shown (personal communication) that one can derive existence of solutions to (54) by considering a related linear RDE. Interestingly, all those solutions have exponential right tail and hence by Proposition 48 none are endogenous.

## 6. Frozen percolation process on infinite binary tree.
A different setting where a particular "max-type" RDE plays the crucial role is the *frozen percolation* process on the infinite binary tree, studied in [5]. Let $\mathbb{T}_3 = (\mathcal{V}, \mathcal{E})$ be the infinite binary tree, where each vertex has degree 3; $\mathcal{V}$ is the set of vertices and $\mathcal{E}$ is the set of undirected edges. Let $(U_e)_{e \in \mathcal{E}}$ be independent edge weights with Uniform$(0, 1)$ distribution. Consider a collection of random subsets $\mathcal{A}_t \subseteq \mathcal{E}$ for $0 \le t \le 1$, whose evolution is described informally by:

$(*)$ $\mathcal{A}_0$ is empty; for each $e \in \mathcal{E}$, at time $t = U_e$ set $\mathcal{A}_t = \mathcal{A}_{t-} \cup \{e\}$ if each end-vertex of $e$ is in a finite cluster of $\mathcal{A}_{t-}$; otherwise set $\mathcal{A}_t = \mathcal{A}_{t-}$.

(A *cluster* is formally a connected component of edges, but we also consider it as the induced set of vertices.) For comparison purposes, a more familiar process is $\mathcal{B}_t := \{e \in \mathcal{E} | U_e \le t\}$, for $0 \le t \le 1$; which gives the standard percolation process on $\mathbb{T}_3$ [35]. It is elementary that the clusters of $\mathcal{B}_t$ can be described in terms of the Galton–Watson branching process and that infinite clusters exist for $t > \frac{1}{2}$ but not for $t \le \frac{1}{2}$. The evolution of the process $(\mathcal{B}_t, 0 \le t \le 1)$ can be described informally by:

$$\text{for each } e \in \mathcal{E}, \text{ at time } t = U_e \text{ set } \mathcal{B}_t = \mathcal{B}_{t-} \cup \{e\}.$$

We notice that any process satisfying $(*)$ must have $\mathcal{A}_t = \mathcal{B}_t$ for $t \le \frac{1}{2}$ but $\mathcal{A}_t \subseteq \mathcal{B}_t$ for $t > \frac{1}{2}$. Qualitatively, in the process $(\mathcal{A}_t)$ the clusters may grow to infinite size but, at the instant of becoming infinite, they are "frozen" in the sense that no extra edge may be connected to an infinite cluster. The final set $\mathcal{A}_1$ will be a forest on $\mathbb{T}_3$ with both infinite and finite clusters, such that no two finite clusters are separated by a single edge.

Following [5] we call this process the *frozen percolation* process on $\mathbb{T}_3$. Although this process is intuitively quite natural, rigorously speaking it is not clear that it exists or that $(*)$ does specify a unique process. In fact Itai Benjamini and Oded Schramm (personal communication) have an argument that such a process does not exist on the $\mathbb{Z}^2$-lattice with its natural invariance property. But for the infinite binary tree case [5] gives a rigorous construction of a process satisfying $(*)$, which can be summarized as follows (Theorem 1 and Proposition 2 of [5]).

THEOREM 52. *There exists a joint law for $(\mathcal{A}_t, 0 \le t \le 1)$ and $(U_e, e \in \mathcal{E})$ such that $(*)$ holds and the joint law is invariant under the automorphisms of $\mathbb{T}_3$. Furthermore for a prescribed edge $e$ and vertex $v$ of $\mathbb{T}_3$, and fixed $t$ in $(\frac{1}{2}, 1)$, the following are true:*

(a) $P(cluster\ containing\ e\ becomes\ infinite\ in\ [t, t+dt]) = \frac{1}{4t^4}\, dt$,

(b) $P(cluster\ containing\ v\ becomes\ infinite\ in\ [t, t+dt]) = \frac{3}{8t^4}\, dt$,

(c) $P(e\ in\ some\ infinite\ cluster\ of\ \mathcal{A}_1) = \frac{7}{12}$, $P(e\ in\ some\ finite\ cluster\ of\ \mathcal{A}_1) = \frac{1}{16}$, $P(e \notin \mathcal{A}_1) = \frac{17}{48}$.

(d) $P(v\ in\ some\ infinite\ cluster\ of\ \mathcal{A}_1) = \frac{7}{8}$, $P(v\ in\ some\ finite\ cluster\ of\ \mathcal{A}_1) = \frac{7}{64}$, $P(v \notin \mathcal{A}_1) = \frac{1}{64}$.

6.1. $540°$ *arguments.* The phrase *circular argument* has negative connotations, but we will describe what we term a $540°$ (i.e., one and a half circles) argument. In summary, the three half-circles are:

- Suppose a process with desired qualitative properties exists. Do heuristic calculations leading to an RDE.
- Solve the RDE. Use the associated RTP to make a rigorous construction of a process.
- Repeat original calculations rigorously.

In the next three sections we outline how this argument is used to prove Theorem 52. A similar $540°$ argument in a more sophisticated setting underlies the mean-field minimal matching example of Section 7.3.

6.2. *Stage* 1. Suppose that the frozen percolation process exists on $\mathbb{T}_3$ and has the natural invariance and independence properties. Define a modified tree called the *planted* binary tree, written $\widetilde{\mathbb{T}}_3 = (\widetilde{\mathcal{V}}, \widetilde{\mathcal{E}})$, where one distinguished vertex (we call it the *root*) has degree 1 and the other vertices have degree 3. Write $\tilde{e}$ for the edge at the root. Clearly $\widetilde{\mathbb{T}}_3$ is isomorphic to the subtree of $\mathbb{T}_3$ which can be obtained by first making some vertex the "root" and then removing two edges coming out of the root and their induced subtrees. Given independent Uniform$(0, 1)$ variables, say $(U_e)_{e \in \widetilde{\mathcal{E}}}$, as the edge weights on $\widetilde{\mathbb{T}}_3$, we suppose we can define a frozen percolation process on this modified tree. Let $Y$ be the time at which the component containing the edge $\tilde{e}$ becomes infinite, with $Y = \infty$ if never. Let $e_1$ and $e_2$ be the two edges which are coming out of the edge $\tilde{e}$; write the corresponding induced planted subtrees as $\widetilde{\mathbb{T}}_{3_1}$ and $\widetilde{\mathbb{T}}_{3_2}$. Let $Y_1$ and $Y_2$ be the respective times for the edges $e_1$ and $e_2$ to join an infinite cluster in $\widetilde{\mathbb{T}}_{3_1}$ or $\widetilde{\mathbb{T}}_{3_2}$. Now consider $\widetilde{\mathbb{T}}_3$. If $U_{\tilde{e}} < \min(Y_1, Y_2)$, then the edge $\tilde{e}$ joins an infinite component; otherwise it never enters the process. Thus one can write

$$(65) \qquad Y = \Phi(\min(Y_1, Y_2), U_{\tilde{e}}),$$

where $\Phi : I \times [0, 1] \to I$, with $I = [\frac{1}{2}, 1] \cup \{\infty\}$ is defined as

$$(66) \qquad \Phi(x, u) = \begin{cases} x, & \text{if } x > u, \\ \infty, & \text{if } x \leq u. \end{cases}$$

Observe that the subtrees $\widetilde{\mathbb{T}}_{3_1}, \widetilde{\mathbb{T}}_{3_2}$ are isomorphic to $\widetilde{\mathbb{T}}_3$, and so $Y_1$ and $Y_2$ are independent and distributed as $Y$. And so the law of $Y$ on the set $I$ satisfies the RDE

$$(67) \qquad Y = \Phi(\min(Y_1, Y_2), U),$$

where $U \overset{d}{=} \mathrm{Uniform}(0,1)$ and $Y_1, Y_2$ are i.i.d. and have the same law as $Y$. Fortunately this RDE is easy to solve.

LEMMA 53 ([5]). *A probability law $\mu$ on $I$ satisfies the RDE* (67) *if and only if for some $x_0 \in [\frac{1}{2}, 1]$*

$$(68) \qquad \mu(dx) = \frac{1}{2x^2}\, dx, \qquad \frac{1}{2} < x \leq x_0;\ \mu(\infty) = \frac{1}{2x_0}.$$

PROOF. Suppose that a probability law $\mu$ on $I$ is a solution of (67) with distribution function $F$. Then from the definition of $\Phi$

$$F(x) = P(U < \min(Y_1, Y_2) \leq x), \qquad \tfrac{1}{2} \leq x \leq 1.$$

Supposing $F$ has a density $F'$ on $[\frac{1}{2}, 1]$ (which can be proved by a more careful rephrasing of the argument); we get

$$F'(x) = 2x(1 - F(x))F'(x), \qquad \tfrac{1}{2} \leq x \leq 1,$$

and hence it follows that

$$(69) \qquad F(x) = 1 - \frac{1}{2x} \qquad \text{on } [\tfrac{1}{2}, 1] \cap \mathrm{support}(\mu).$$

Since the function $x \mapsto 1 - \frac{1}{2x}$ is strictly increasing, identity (69) can only happen when $\mathrm{support}(\mu) = [\frac{1}{2}, x_0]$ for some $\frac{1}{2} < x_0 \leq 1$.

Conversely it is easy to see that such a probability law on $I$ satisfies the RDE (67). $\square$

From the definition of $Y$ in terms of frozen percolation on the planted binary tree, we expect the support of its distribution to be all of $I = [\frac{1}{2}, 1] \cup \{\infty\}$, and so we choose the particular solution (68) with $x_0 = 1$, that is, the distribution $\nu$ defined by

$$(70) \qquad \nu(dy) = \frac{1}{2y^2}\, dy, \qquad \frac{1}{2} \leq y \leq 1,\ \nu(\infty) = \frac{1}{2},$$

or equivalently

$$(71) \qquad \nu((y, \infty]) = \frac{1}{2y}, \qquad \frac{1}{2} \leq y \leq 1.$$

Continuing to argue heuristically, we can now do the calculations recorded in Theorem 52: we give the argument for (a), and the other parts are similar.

Write $e_1, e_2, e_3, e_4$ for the edges adjacent to the edge $e$, and $\widetilde{\mathbb{T}}_{3_1}, \widetilde{\mathbb{T}}_{3_2}, \widetilde{\mathbb{T}}_{3_3}, \widetilde{\mathbb{T}}_{3_4}$ for the corresponding planted binary trees which are all isomorphic to $\widetilde{\mathbb{T}}_3$, and further let $Y_1, Y_2, Y_3, Y_4$ be the times at which the respective edges enter an infinite cluster of the frozen percolation processes restricted to the subtrees. Writing $Z$ for the time taken for the edge $e$ to enter an infinite cluster (note $Z = \infty$ if never), then

$$(72) \qquad Z = \begin{cases} \min\limits_{1 \le i \le 4} Y_i, & \text{if } U_e < \min\limits_{1 \le i \le 4} Y_i, \\ \infty, & \text{otherwise.} \end{cases}$$

Thus the density $f_Z$ of $Z$ on $[\frac{1}{2}, 1]$ in terms of the law $\nu$ of $Y$ can be written as

$$f_Z(x) = x \times 4 \frac{d\nu}{dx} \times \nu^3((x, \infty)) = 4x \times \frac{1}{2x^2} \times \left(\frac{1}{2x}\right)^3 = \frac{1}{4x^4},$$

as asserted in part (a).

6.3. *Stage* 2. We now start a rigorous construction based on knowing that the distribution $\nu$ at (70) solves the RDE (67). For each edge $e \in \mathcal{E}$ let $\vec{e}$ and $\overleftarrow{e}$ be the two directed edges defined by $e$, and let $\vec{\mathcal{E}}$ be the set of all directed edges. Now for the directed edges we have a natural language of family relationship: the edge $\vec{e} = (v, w)$ has two *children* of the form $(w, x_1)$ and $(w, x_2)$. It is not hard to use the Kolmogorov consistency theorem and the fact that $\nu$ solves (67) to show

LEMMA 54 ([5]). *There exists a joint law for* $((U_{\vec{e}}, Y_{\vec{e}}), \vec{e} \in \vec{\mathcal{E}})$ *which is invariant under automorphisms of* $\mathbb{T}_3$ *and such that for each* $\vec{e} \in \vec{\mathcal{E}}$

$$(73) \qquad \begin{aligned} &Y_{\vec{e}} \text{ has law } \nu, \\ &Y_{\vec{e}} = \Phi(\min(Y_{\vec{e}_1}, Y_{\vec{e}_2}), U_{\vec{e}}) \qquad a.s., \end{aligned}$$

*where* $\vec{e}_1$ *and* $\vec{e}_2$ *are children of* $\vec{e}$, *and where for each* $e \in \mathcal{E}$, $U_{\vec{e}} = U_{\overleftarrow{e}} = U_e$.

Now we can outline the rigorous construction of the frozen percolation process. Essentially, one takes the heuristically obvious property (72) as a definition. In more detail, for an undirected edge $e$, define $\partial_e$ as the set of four directed edges adjacent to $e$ and directed away from it. Define

$$(74) \qquad \mathcal{A}_1 := \{e \in \mathcal{E} | U_e < \min(Y_{e'} : e' \in \partial_e)\}.$$

Finally for $0 \le t < 1$ define

$$(75) \qquad \mathcal{A}_t := \{e \in \mathcal{A}_1 | U_e \le t\}.$$

It is now clear that $(\mathcal{A}_t)$ inherits from $(Y_{\vec{e}})$ the automorphism-invariance property, as well as the property that the only possible time at which an edge $e$ can join the process is at time $U_e$. To check that $(\mathcal{A}_t)$ actually satisfies (∗) and so it is a frozen percolation process is somewhat more elaborate, and this part of the proof [5] of existence is omitted here.

6.4. *Stage* 3. We can now repeat the argument leading to (72) in terms of the explicit construction above and its modification on planted binary trees. This leads to part (a) as shown in Stage 1, and the other parts are similar.

6.5. *The endogenous property.* If the RDE were nonendogenous, then the frozen percolation process would have a kind of "spatial chaos" property, that the behavior near the root was affected by the behavior at infinity. For several years we conjectured in seminar talks that the RDE is nondogenous, but recently proved the opposite.

THEOREM 55 ([16]). *The invariant RTP associated with* (67), (70) *is endogenous.*

**7. Combinatorial optimization within the mean-field model of distance.** In problems involving $n$ random points in $d$-dimensional space, explicit calculations for $n \to \infty$ asymptotics are often complicated by the obvious fact that the $\binom{n}{2}$ inter-point distances are dependent r.v.'s. One can make a less realistic but more tractable model by eliminating the ambient $d$-dimensional space and instead assuming that the $\binom{n}{2}$ inter-point distances are *independent* r.v.'s. This is the *mean-field model of distance*. Specifically, assume inter-point distances have exponential distribution with mean $n$, so nearest-neighbor distances are order 1. This model mimics true inter-point distances in $d = 1$ dimension; other distributions can be used to mimic other $d$ without changing essential aspects of what follows.

This model, and study of the *minimal spanning tree* and *minimal matching* problems within it, are surveyed in some detail in Sections 4 and 5 of [11]. Here we emphasize a different example, in Section 7.2, and only briefly record the RDEs arising in the minimal matching, traveling salesman and variant problems (Sections 7.3–7.6).

7.1. *The PWIT approximation.* In the mean-field model above, the key feature is that there is an $n \to \infty$ "local weak limit" structure called the PWIT (*Poisson weighted infinite tree*), which describes the geometry of the space as seen from a fixed reference point. In brief (see [11] for more details) consider a Poisson point process

$$(76) \qquad\qquad 0 < \xi_1 < \xi_2 < \xi_3 < \cdots$$

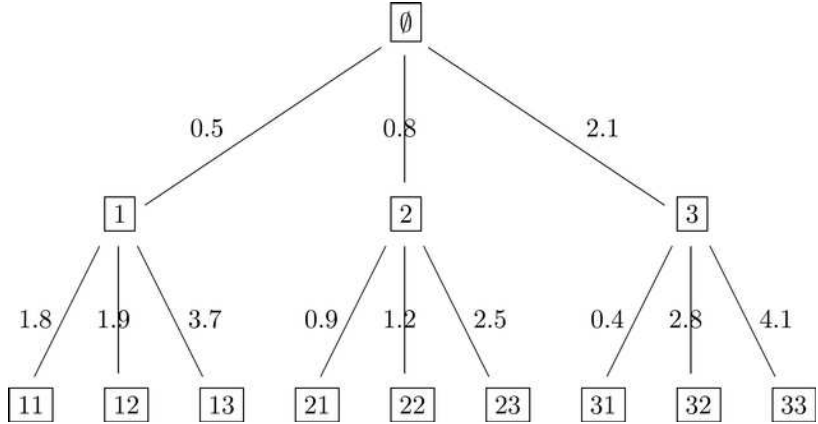

Fig. 4. *Part of a realization of the PWIT that shows just the first three children of each vertex. The length is written next to each edge e.*

of rate 1 on $(0, \infty)$. Take a root vertex $\varnothing$. Let this root have an infinite number of children $1, 2, 3, \ldots$, the edge-lengths to these children being distributed as the Poisson process $(\xi_i, i \geq 1)$ at (76). Repeat recursively; each vertex $\mathbf{i}$ has an infinite number of children $(\mathbf{i}j, \; j \geq 1)$ and the edge-lengths $\xi_{\mathbf{i}j}, \; j \geq 1$ are distributed as the Poisson process (76), independent of other such Poisson processes. See Figure 4.

7.2. *Critical point for minimal subtrees.* Consider the mean-field model on $n$ points as the complete graph $K_n$ on $n$ vertices, and write $\xi_e$ for the length of edge $e$. For a subtree $\mathbf{t}$, that is, a tree whose vertices are some subset of the $n$ vertices, write

$$|\mathbf{t}| = \text{ number of edges in } \mathbf{t},$$

$$L(\mathbf{t}) = \sum_{e \in \mathbf{t}} \xi_e = \text{ total length of } \mathbf{t},$$

$$a(\mathbf{t}) = L(\mathbf{t})/|\mathbf{t}| = \text{ average edge-length of } \mathbf{t}.$$

A well-known result [32] on minimal spanning trees says that, if we insist on $|\mathbf{t}| = n - 1$, then the smallest we can make $a(\mathbf{t})$ is about $\zeta(3) := \sum_i i^{-3}$. If we fix $0 < \varepsilon < 1$ and consider subtrees with around $\varepsilon n$ edges, then we guess that the smallest value of $a(\mathbf{t})$ should be around $\delta(\varepsilon)$ as $n \to \infty$, for some deterministic function $\delta(\varepsilon)$. It is not hard to show that $\delta(\varepsilon) > 0$ for large $\varepsilon$ while $\delta(\varepsilon) = 0$ for small $\varepsilon$. So there must be a *critical point* at which $\delta(\cdot)$ leaves 0; this is analogous to critical points in percolation theory. What is interesting is that the critical point is determined by an RDE. It is convenient to turn the problem around and study the *maximum* size of a subtree $\mathbf{t}$ subject to the constraint that $a(\mathbf{t}) \leq c$. For fixed $0 < c < \infty$ consider the RDE on $S = [0, \infty)$

$$Y \overset{d}{=} \sum_{i=1}^{\infty} (c - \xi_i + Y_i)^+; \qquad (\xi_i) \text{ a Poisson rate 1 point process on } (0, \infty).$$
(77)

PROPOSITION 56 ([4]).  *Define* $M(n, c) = \max\{|\mathbf{t}| : \mathbf{t}$ *a subtree of* $K_n, a(\mathbf{t}) \leq c\}$. *Then there exists a critical point* $c(0) \in [e^{-2}, e^{-1}]$ *such that*

(78)  *if* $c < c(0)$        *then* $n^{-1}M(n, c) \overset{d}{\to} 0,$

(79)  *if* $c > c(0)$        *then* $\exists \eta(c) > 0$ *such that* $P(n^{-1}M(n, c) \geq \eta(c)) \to 1,$

*and*

$$c(0) = \inf\{c \,|\, \text{RDE (77) has no solution on } [0, \infty)\}.$$

The conceptual point to emphasize is that, by analogy with Example 2, we are studying an "average" by studying whether an associated "compensated sum" is finite or infinite.

First we explain how the RDE (77) arises. In the PWIT define, for integers $h \geq 0$,

$$Y^{(h)} = \sup\{c|\mathbf{t}| - L(\mathbf{t})|\text{root} \in \mathbf{t}, \text{height}(\mathbf{t}) \leq h\},$$

where the *sup* is over subtrees $\mathbf{t}$ of the PWIT, and where height($\mathbf{t}$) denotes the maximum number of edges in a path in $\mathbf{t}$ from the root. To obtain the maximizing $\mathbf{t}$ one simply considers in turn each child $i$ of the root and considers whether one gets a positive contribution by including child $i$ in $\mathbf{t}$. The contribution equals

$$c - \xi_i + Y_i^{(h-1)}$$

where $\xi_i$ is the length of edge from root to child $i$, and $Y_i^{(h-1)}$ is a sup over subtrees $\mathbf{t}_i$ of child $i$:

$$Y_i^{(h-1)} = \sup\{c|\mathbf{t}_i| - L(\mathbf{t}_i)|\text{height}(\mathbf{t}_i) \leq h-1\}.$$

So

(80) $$Y^{(h)} = \sum_{i=1}^{\infty} (c - \xi_i + Y_i^{(h-1)})^+$$

where the $(Y_i^{(h-1)})$ are, by the recursive structure of the PWIT, independent copies of $Y^{(h-1)}$. Writing $T_c$ for the map on distributions associated with the RDE (77), the last equality says

$$T_c^h(\delta_0) = \text{dist}(Y^{(h)}).$$

Lemma 15 then implies that for $c < c(0)$,

$$T_c^h(\delta_0) \xrightarrow{d} \mu_c \qquad \text{as } h \to \infty,$$

where $\mu_c$, supported on $[0, \infty)$, is the lower invariant measure; and that the RTP is endogenous. Indeed, $\mu_c$ is just the distribution of

$$Y^{(\infty)} = \sup\{c|\mathbf{t}| - L(\mathbf{t})|\text{root} \in \mathbf{t}, \ \mathbf{t} \text{ finite}\}.$$

OUTLINE PROOF OF PROPOSITION 56. Fix $c < c(0)$. Roughly, the fact that $Y^{(\infty)}$ is finite implies that there cannot exist large subtrees of the PWIT with average edge-length much greater than $c$; the fact that the PWIT represents the local structure of $K_n$ for large $n$ implies that the same should hold for $K_n$; this is the lower bound (78) of the proposition. Let us amplify this argument into four steps. Fix an integer $m$.

*Step* 1. The connection between $K_n$ (the complete graph on $n$ vertices with random edge-lengths) and the PWIT is provided by *local weak convergence* of the former to the latter—see [11], Theorem 3 for formalization. A soft consequence of local weak convergence is

(81)     $$\lim_n \text{dist} \max\{c|\mathbf{t}| - L(\mathbf{t})|\mathbf{t} \subset K_n, \text{root} \in \mathbf{t}, \ |\mathbf{t}| \leq 3m\}$$
        is stochastically smaller than $T_c^{3m}\delta_0$.

Indeed, we would have asymptotic equality if we required only that $\mathbf{t}$ have depth $\leq 3m$; but we make a stronger restriction.

*Step* 2. The quantity above can be used to bound the chance of the event: there exists a small tree $\mathbf{t}$ containing the root and with $c|\mathbf{t}| - L(\mathbf{t}) \geq x$. Consider the mean number of vertices $v$ satisfying that event (with $v$ in place of root) and apply Markov's inequality to deduce the following:

$P(\exists \text{ at least } \delta n \text{ vertices } v \text{ s.t. } \exists \mathbf{t} \ni v \text{ with } c|\mathbf{t}| - L(\mathbf{t}) \geq x \text{ and } |\mathbf{t}| \leq 3m)$
        is asymptotically $\leq \delta^{-1}\mu_c[x, \infty]$.

*Step* 3. Any tree with at least $3m$ edges can be split into edge-disjoint subtrees, each having between $m$ and $3m$ edges.

*Step* 4. We are assuming $c < c(0)$, so interpose $c < c' < c^* < c(0)$. If the desired conclusion (78) were false, there would be some "big" tree $\mathbf{t}$ with $|\mathbf{t}| \geq \varepsilon n$ and $a(\mathbf{t}) \leq c$. Use step 3 to split into "small" subtrees; a deterministic averaging argument shows there would exist at least $\delta n$ vertices $v$ in small trees $\mathbf{t}'$ with $a(\mathbf{t}') \leq c'$. Here $\delta$ depends on $\varepsilon, c, c'$ but not on $n, m$. These trees now satisfy $c^*|\mathbf{t}| - L(\mathbf{t}) \geq m(c^* - c')$. Applying step 2, the chance of this many small trees is at most $\delta^{-1}\mu_{c^*}[m(c^* - c'), \infty]$. Since $m$ is arbitrary and $\mu_{c^*}(\infty) = 0$ we get the lower bound (78). We are done.

The argument for the upper bound (79) is parallel. For $c > c(0)$ each $v$ has nonvanishing chance of being in some large finite tree $\mathbf{t}$ with $a(\mathbf{t}) \geq \varepsilon(c)$, and one can patch together these trees to get an $\Omega(n0$-size tree $\mathbf{t}$ with $a(\mathbf{t})) \geq \varepsilon(c)$. $\square$

Incidentally, the lower bound $c(0) \geq e^{-2}$ stated in Proposition 56 arises from the first moment method, and the upper bound $c(0) \leq e^{-1}$ comes from considering *paths* as a special case of trees. Moreover, numerically $c(0) \approx 0.263$. See [4] for details.

7.3. *Minimal matching.* Consider again the mean-field model of distance, that is, the complete graph $K_n$ with random edge-lengths with exponential (mean $n$) distribution. Take $n$ even and consider a (complete) matching, that is, a collection of $n/2$ vertex-disjoint edges. Define

$$M_n = \text{minimum total length of a complete matching.}$$

This problem is often studied in the bipartite case ([61], Chapter 4) but the two versions turn out to be equivalent in our asymptotic setting. The following limit behavior was argued nonrigorously in [53] and proved (in the bipartite setting) in [1, 6]. There are fascinating recent proofs [48, 55] of an underlying exact formula for $\mathbb{E}M_n$ in the bipartite, exponential distribution setting, but it seems unlikely that the applicability of exact methods extends far into the broad realm of problems amenable to asymptotic study.

THEOREM 57. $\frac{2}{n}\mathbb{E}M_n \to \pi^2/6$.

The technically difficult proof is outlined in moderate detail in Section 5 of [11]. Here we emphasize only the underlying RDE, and some analogous RDEs arising in analogous problems.

The central idea is that, since the PWIT originates as a "local weak limit" of $K_n$, one can relate matchings on $K_n$ to matchings on the PWIT. The technically hard, though noncomputational, part of the proof is to show that the limit $\lim_n \frac{2}{n}\mathbb{E}M_n$ must equal

$$\text{inf}\{\mathbb{E}(\text{typical edge-length of } \mathcal{M}) : \mathcal{M} \text{ an invariant matching on the PWIT}\}. \tag{82}$$

Here *invariant* means, intuitively, that in defining the matching $\mathcal{M}$ on the PWIT, the root $\varnothing$ must play no special role. Now one can see how to construct the optimal matching $\mathcal{M}_{\text{opt}}$ on the PWIT by reusing two ideas we have seen earlier in this survey. First, we use the 540° *argument* from Section 6.1: start with heuristically defined quantities, obtain an RDE and use its solution as a basis for rigorous construction. Second, we use the idea from Section 4.6 of seeking a recursion for a quantity defined as a *difference*.

Write **T** for the PWIT. Consider the definition, analogous to (50),

(83)        $X_\varnothing =$ length of optimal matching on **T**
                     $-$ length of optimal matching on $\mathbf{T} \setminus \{\varnothing\}$.

Here we mean *total* length, so we get $\infty - \infty$, and so this makes no sense rigorously. But pretend it does make sense. Then for each child $j$ of the root we can define $X_j$ similarly in terms of the subtree $\mathbf{T}^j$ rooted at $j$:

        $X_j =$ length of optimal matching on $\mathbf{T}^j$
                $-$ length of optimal matching on $\mathbf{T}^j \setminus \{j\}$.

One can now argue, analogously to (51),

(84)                   $$X_\varnothing = \min_{1 \le j < \infty} (\xi_j - X_j),$$

(85) root is matched to the vertex $\arg\min_j (\xi_j - X_j)$ in the optimal matching.

Recall $(\xi_i)$ is the Poisson process (76). This motivates us to consider the RDE

(86)                 $$X \stackrel{d}{=} \min_{1 \le i < \infty} (\xi_i - X_i) \qquad (S = \mathbb{R}).$$

Luckily, this turns out to be easy to solve.

LEMMA 58 ([6], Lemma 5).   *The unique solution of* (86) *is the logistic distribution*

(87)          $$P(X \le x) = 1/(1 + e^{-x}), \qquad -\infty < x < \infty,$$

*or equivalently the density function*

              $$f(x) = (e^{x/2} + e^{-x/2})^{-2}, \qquad -\infty < x < \infty.$$

Implementing the 540° argument, we will use the logistic solution or the RDE to construct a random matching on the PWIT. Each edge $e$ in the edge-set $\mathbb{E}$ of **T** corresponds to two directed edges $\vec{e}, \overleftarrow{e}$: write $\vec{\mathbb{E}}$ for the set of directed edges and write $\xi(\vec{e}) = \xi(\overleftarrow{e}) = \xi(e)$ for the edge-length. For a directed edge $(v, w)$ we can call the directed edges $\{(w, x)| \ x \ne v\}$ its *children*. The Kolmogorov consistency theorem and the logistic solution of (86) imply (cf. Lemma 54)

LEMMA 59.   *Jointly with the edge-lengths* $(\xi(e), e \in \mathbb{E})$ *of the PWIT we can construct* $\{X(\vec{e}), \vec{e} \in \vec{\mathbb{E}}\}$ *such that:*

   (i)  *each* $X(\vec{e})$ *has the logistic distribution,*

(ii) *for each $\vec{e}$, with children $\vec{e}_1, \vec{e}_2, \dots$ say,*

$$X(\vec{e}) = \min_{1 \le j < \infty}(\xi(e_j) - X(\vec{e}_j)). \tag{88}$$

Theorem 61 will show that $X(v, v')$ depends only on the edge-lengths within the subtree rooted at $v'$. Guided by the heuristic (85), for each vertex $v$ define

$$v^* = \arg\min_{v' \sim v}(\xi(v, v') - X(v, v')). \tag{89}$$

In view of (82) an outline proof of Theorem 57 can be completed by proving

PROPOSITION 60 ([6], Lemma 16, Propositions 17 and 18).  (a) *The set of edges $(v, v^*)$ forms a matching $\mathcal{M}_{\mathrm{opt}}$ on the PWIT.*

*Write $\vec{\mathcal{M}}_{\mathrm{opt}}(\varnothing)$ for the vertex to which the root $\varnothing$ is matched in $\mathcal{M}_{\mathrm{opt}}$, so that the mean edge-length in $\mathcal{M}_{\mathrm{opt}}$ can be written as $\mathbb{E}\xi(\varnothing, \vec{\mathcal{M}}_{\mathrm{opt}}(\varnothing))$. Then:*

(b) $\mathbb{E}\xi(\varnothing, \vec{\mathcal{M}}_{\mathrm{opt}}(\varnothing)) = \pi^2/6$.

(c) *For any invariant matching $\mathcal{M}$,*

$$\mathbb{E}\xi(\varnothing, \vec{\mathcal{M}}(\varnothing)) - \mathbb{E}\xi(\varnothing, \vec{\mathcal{M}}_{\mathrm{opt}}(\varnothing)) \ge 0. \tag{90}$$

Let us indicate only the proofs of (a) and (b). For (a) we need only show that $(v^*)^* = v$. Note first

$$\xi(v, v^*) - X(v, v^*) < \min_{y \ne v}(\xi(v, y) - X(v, y)) \qquad \text{by definition of } v^*$$

$$= X(v^*, v) \qquad\qquad\qquad \text{by recursion (88)}$$

or equivalently

$$\xi(v, v^*) < X(v, v^*) + X(v^*, v). \tag{91}$$

And if $z \ne v^*$ is another neighbor of $v$, then

$$\xi(v, z) - X(v, z) > \min_{y \ne v}(\xi(v, y) - X(v, y))$$

$$= X(z, v)$$

or equivalently

$$\xi(v, z) > X(v, z) + X(z, v).$$

We conclude that $v^*$ is the *unique* neighbor of $v$ satisfying (91). But the right-hand side of (91) is symmetric, so applying this conclusion to $v^*$ shows $(v^*)^* = v$.

To prove (b) we calculate the mean length of the edge at the root. In order for this length to be $x$, there must be an edge of length $x$ from the root to some vertex $j$, and also (91) we must have $x < X(\text{root}, j) + X(j, \text{root})$. But these are distributed as independent logistics, say $X_1$ and $X_2$, and so

$$(92) \quad \mathbb{E}\xi(\varnothing, \vec{\mathcal{M}}_{\text{opt}}(\varnothing)) = \int_0^\infty x \, dx \; P(x < X_1 + X_2)$$

$$= \tfrac{1}{2}\mathbb{E}((X_1 + X_2)^+)^2 \qquad \text{by a general formula}$$

$$= \tfrac{1}{4}\mathbb{E}(X_1 + X_2)^2 \qquad \text{by symmetry}$$

$$= \tfrac{1}{2}\mathbb{E}X_1^2$$

$$(93) \qquad = \pi^2/6,$$

the last step using a standard fact that the logistic distribution has variance $\pi^2/3$.

THEOREM 61 ([15]).  *The invariant RTP associated with the RDE* (86) *is endogenous.*

The significance of this result is pointed out in Section 7.5. The proof involves somewhat intricate analytic study of the iterates $T^{(2)n}(\mu \times \mu)$ to verify Theorem 11(c). We remark that we have not succeeded in using contraction methods to prove Theorem 61. Indeed the operator $T$ associated with the RDE (86) is not a strict contraction. To see this, it is easy to check that $T$ is well defined on the subspace $\mathcal{P}_1$ of distributions with finite mean. Moreover, if $X$ has logistic distribution, then the distribution of $c + X$ is a fixed point for $T^2$ for any $c \in \mathbb{R}$. Hence $T$ cannot be a strict contraction on the whole of $\mathcal{P}_1$. This shows that the logistic solution of (86) does not have full "domain of attraction," but the question of determining the domain of attraction remains open.

OPEN PROBLEM 62.  Find the subset $\mathcal{D} \subset \mathcal{P}_1$ such that $T^n(\nu) \xrightarrow{d} \mu$ as $n \to \infty$ if and only if $\nu \in \mathcal{D}$, where $\mu$ is the logistic distribution.

REMARK.  The way we started the heuristic argument at (83) may seem implausible, because one might expect analogous quantities in finite-$n$ setting to have spreads which increase to infinity with $n$. But a better analogy is with the position $R_n$ of the rightmost point in BRW; its spread (Lemma 43) stays bounded with $n$.

7.4. *TSP and other matching problems.* Here we follow Sections 6.1 and 6.2 of [6], which showed how earlier work [51, 52] fits into the current framework.

As suggested at the start of Section 7, one can define a mean-field model of distance with any real *pseudo-dimension* $0 < d < \infty$ to mimic distances between random points in $d$ dimensions. Precisely, take the complete graph $K_n$ on $n$ vertices, and let the i.i.d. edge-lengths have distribution $n^{1/d}L$ where

$$P(L \le x) \sim x^d/d \qquad \text{as } x \downarrow 0.$$

This scaling keeps nearest-neighbor distances as order 1. And in the local approximation of this $K_n$ by the PWIT, we simply change the distribution of edge-lengths at a vertex in the PWIT: the homogeneous Poisson process of rate 1 is replaced by an inhomogeneous Poisson process $0 < \xi_1 < \xi_2 < \cdots$ of rate $x^{d-1}$; in other words

$$\mathbb{E}(\text{number of } i \text{ with } \xi_i \le x) = x^d/d.$$

For minimum matching in pseudo-dimension $d$, it is remarkable that (heuristically, at least) the previous analysis is in principle unchanged. Theorem 57 becomes [cf. (92)]

$$\lim_n \frac{2}{n}\mathbb{E}M_n = \int_0^\infty x^d P(X_1 + X_2 > x)\,dx$$

where $X_1$ and $X_2$ are independent copies of the solution of the RDE

$$(94) \qquad X \stackrel{d}{=} \min_{1 \le i < \infty} (\xi_i - X_i) \qquad (S = \mathbb{R}).$$

Thus the abstract structure of the limit theorem is unchanged in pseudo-dimension $d$. But for $d \ne 1$ there is no known explicit solution of (94); and while numerical methods indicate that there is indeed a unique solution, rigorous proof remains elusive.

OPEN PROBLEM 63. Prove that for each real $0 < d < \infty$ there is a unique solution to the RDE (94), and that the associated invariant RTP is endogenous.

Similarly, in the TSP (traveling salesman problem) in pseudo-dimension $d$, a variant of the argument leading to recursion (86) leads us to the recursion

$$(95) \qquad X \stackrel{d}{=} \min_{1 \le i < \infty}^{[2]} (\xi_i - X_i) \qquad (S = \mathbb{R}).$$

Here $\min^{[2]}$ denotes the second minimum, and the analog of Theorem 57 is that the length $S_n$ of the optimal TSP satisfies

$$\lim_n \frac{1}{n}\mathbb{E}S_n = \int_0^\infty x^d P(X_1 + X_2 > x)\,dx$$

where $X_1$ and $X_2$ are independent copies of the solution of the RDE (95). Numerically the limit is about 2.04 for $d = 1$. Again numerical methods indicate that there is a unique solution for all $d$, but no rigorous proof is known even for $d = 1$.

OPEN PROBLEM 64. Prove that for each real $0 < d < \infty$ (or at least for $d = 1$) there is a unique solution to the RDE (95), and that the associated invariant RTP is endogenous.

Instead of studying minimal matchings one could study Gibbs distributions on matchings; this leads to a different RDE ([6], (46), and [62])

$$1/X \overset{d}{=} \sum_{i=1}^{\infty} e^{-\theta \xi_i} X_i \qquad (S = \mathbb{R}^+)$$

which is somewhat in the spirit of the linear case.

7.5. *The cavity method.* The nonrigorous *cavity method* was developed in statistical physics in the 1980s; see [54] for a recent survey. Though typically applied to examples such as ground states of disordered Ising models, it can also be applied to the kind of "mean-field combinatorial optimization" examples of the last two sections. It turns out that the methodology used in [6] to make a rigorous proof of the mean-field matching limit serves to provide a general methodology for seeking rigorous proofs paralleling the cavity method in a variety of contexts. This is a broad and somewhat complicated topic, and the time is not ripe for a definitive survey, but it seems worthwhile to outline the ingredients of the rigorous methodology, pointing out where RDEs and endogeny arise.

Start with a combinatorial optimization problem over some size-$n$ random structure.

- Formulate a "size-$\infty$" random structure, the $n \to \infty$ limit in the sense of local weak convergence.
- Formulate a corresponding combinatorial optimization problem on the size-$\infty$ structure.
- Heuristically define relevant quantities on the size-$\infty$ structure via additive renormalization [cf. (83)].
- If the size-$\infty$ structure is treelike (the only case where one expects exact asymptotic solutions), observe that the relevant quantities satisfy a problem-dependent RDE.
- Solve the RDE. Use the unique solution to find the value of the optimization problem on the size-$\infty$ structure.
- Show that the RTP associated with the solution is endogenous.

- Endogeny shows that the optimal solution is a measurable function of the data, in the infinite-size problem. Since a measurable function is almost continuous, we can pull back to define almost-feasible solutions of the size-$n$ problem with almost the same cost.
- Show that in the size-$n$ problem one can patch an almost-feasible solution into a feasible solution for asymptotically negligible cost.

7.6. *Scaling laws in mean-field combinatorial optimization.* Here we indicate current nonrigorous work on scaling exponents associated with mean field of combinatorial optimization problems. As indicated in the methodology description above, the main requirement for making rigorous proofs would be proofs of uniqueness and endogeny for the RDEs which arise. So our discussion emphasizes the RDEs.

7.6.1. *Near-optimal solutions.* In the context of mean-field minimum matching, compare the optimal matching $\mathcal{M}^{(n)}_{\mathrm{opt}}$ with a near-optimal matching $\mathcal{M}^{(n)}$ by using the two quantities

$$\delta_n(\mathcal{M}^{(n)}) = n^{-1}\mathbb{E}[\text{number edges of } \mathcal{M}^{(n)} \setminus \mathcal{M}^{(n)}_{\mathrm{opt}}],$$

$$\varepsilon_n(\mathcal{M}^{(n)}) = n^{-1}\mathbb{E}[\text{cost}(\mathcal{M}^{(n)}) - \text{cost}(\mathcal{M}^{(n)}_{\mathrm{opt}})].$$

Then define

$$\varepsilon^n_*(\delta) = \min\{\varepsilon_n(\mathcal{M}^{(n)})|\delta_n(\mathcal{M}^{(n)}) \geq \delta\}.$$

We anticipate a limit $\varepsilon(\delta) = \lim_n \varepsilon^n_*(\delta)$, and then can ask whether there is a scaling exponent

$$\varepsilon(\delta) \asymp \delta^\alpha \qquad \text{as } \delta \to 0.$$

Such a scaling exponent provides a measure of how different an almost-optimal solution can be from the optimal solution.

Remarkably, it is not so hard to study this question by an extension of the methods of Section 7.3. It turns out [10] that the key is the extension of the RDE (86) to the following RDE on $S = \mathbb{R}^3$:

$$(96) \qquad \begin{pmatrix} X \\ Y \\ Z \end{pmatrix} = \begin{pmatrix} \min_i(\xi_i - X_i) \\ \min_i(\xi_i - (Z_i + \lambda)\mathbb{1}(i = i^*) - Y_i\mathbb{1}(i \neq i^*)) \\ \min_i(\xi_i - Y_i) \end{pmatrix}$$

where

$$i^* = \arg\min_i(\xi_i - X_i)$$

and where $\lambda > 0$ is a Lagrange multiplier. In terms of the solution of this RDE one can define functions $\varepsilon(\lambda)$ and $\delta(\lambda)$ which then define the limit function $\varepsilon(\delta)$. Numerical study in [10] indicates the scaling exponent $\alpha = 3$ in both minimal matching and TSP problems in the mean-field model.

7.6.2. *TSP percolation function.* In the context of the TSP in the mean-field model of distance, one can study a function $(p(u), 0 < u \le 1)$ analogous to the *percolation function*, defined as follows. Recall $K_n$ is the complete graph with random edge-lengths. Over all cycles $\pi_{n,u}$ in $K_n$ containing $un$ vertices, let $C_{n,u}$ be the minimum average edge-length of $\pi_{n,u}$. We anticipate a limit

$$\lim_n \mathbb{E} C_{n,u} = p(u).$$

It turns out [8] that what is relevant is the following RDE for a distribution $(X, Z)$ on $\mathbb{R}^2$:

$$(97) \quad \begin{pmatrix} X \\ Z \end{pmatrix} \stackrel{d}{=} \begin{pmatrix} \max_i (\lambda - \xi_i + X_i - Z_i^+) \\ \max_i (\lambda - \xi_i + X_i - Z_i^+) + \max_i^{[2]} (\lambda - \xi_i + X_i - Z_i^+) \end{pmatrix}.$$

Here $\lambda > 0$ is again a Lagrange multiplier. In terms of the solution of this RDE one can define functions $p(\lambda)$ and $u(\lambda)$ which then define the limit function $p(u)$. Numerical study in [8] indicates a scaling exponent $p(u) \asymp u^\alpha$ as $u \downarrow 0$ with $\alpha = 3$. Moreover, for both the present "percolation function" setting and the previous "near-optimal solution" setting, one can pose analogous questions involving trees in place of tours, and it turns out [8, 10] that for both these questions the scaling exponent is 2. But at present we have no good conceptual explanation of these fascinating observations.

7.6.3. *First passage percolation.* A somewhat different setup is appropriate for a mean-field model of first passage percolation. Take the 4-regular tree $\mathbb{T}$, with in-degree 2 and out-degree 2 at each vertex; regard this as the mean-field analog of the oriented lattice $\mathbb{Z}^2$. Attach independent exponential(1) random variables $\xi_e$ to the edges of $\mathbb{T}$. We study *flows* $\mathbf{f} = (f(e))$ on $\mathbb{T}$, for which the in-flow equals the out-flow at each vertex, with $0 \le f(e) \le 1$. Associated with an invariant random flow are two numbers

$\partial(\mathbf{f}) = \mathbb{E} f(e)$:         the average density of the flow,

$\tau(\mathbf{f}) = \dfrac{\mathbb{E} f(e) \xi(e)}{\partial(\mathbf{f})}$:     the flow-weighted average edge-traversal time.

We study the function

$\delta^*(\tau) := \sup\{\partial(\mathbf{f}) : \mathbf{f} \text{ an invariant flow with } \tau(\mathbf{f}) = \tau\}, \qquad 0 < \tau < 1.$

We have $\delta^*(\tau) > 0$ iff $\tau > \tau_{\mathrm{FPP}}$, where $\tau_{\mathrm{FPP}}$ is the time constant in first passage percolation on $\mathbb{T}$. As above, to study scaling exponents the key is a certain RDE for $S = \mathbb{R}^+$, which turns out [7] to be

$$(98) \quad Z \stackrel{d}{=} \min\left(Z_2 + \xi_2 - a, Z_3 + \xi_3 - a, \sum_{i=1}^{3}(Z_i + \xi_i - a)\right)$$
$$- \min\left(0, \sum_{i=1,2}(Z_i + \xi_i - a), \sum_{i=1,3}(Z_i + \xi_i - a)\right).$$

Here $a$ is a parameter $\in (\tau_{\mathrm{FPP}}, 1)$. In terms of the solution of this RDE one can define functions $\delta(a)$ and $\tau(a)$ which then determine the function $\delta^*(\tau)$. Numerical study in [7] indicates a scaling exponent 2:

$$\delta^*(\tau) \sim 12.7(\tau - \tau_{\mathrm{FPP}})^2 \quad \text{as } \tau \downarrow \tau_{\mathrm{FPP}}.$$

## 8. Complements.

8.1. *Numerical and Monte Carlo methods.* In the context of studying a fixed point equation $T(\mu) = \mu$ or the bivariate analog in Theorem 11, there are several numerical methods one might try: solving the equation directly or calculating iterates $T^n(\mu_0)$ for some convenient $\mu_0$; discretization or working in a basis expansion. But implementation is highly problem-dependent.

In contrast, given an RDE $X \stackrel{d}{=} g(\xi, X_i, \ i \geq 1)$ the *bootstrap Monte Carlo* method provides a very easy to implement and essentially problem-independent method. Start with a list of $N$ numbers (take $N = 10,000$ say) with empirical distribution $\mu_0$. Regard these as "generation-0" individuals $(X_i^0, 1 \leq i \leq N)$. Then $T(\mu_0)$ can be approximated as the empirical distribution $\mu_1$ of $N$ "generation-1" individuals $(X_i^1, \ 1 \leq i \leq N)$, each obtained independently via the following procedure. Take $\xi$ with the prescribed distribution, take $I_1, I_2, \ldots$ independent uniform on $\{1, 2, \ldots, N\}$ and set

$$X_i^1 = g(\xi, X_{I_1}^0, X_{I_2}^0, \ldots).$$

Repeating for, say, 20 generations lets one see whether $T^n(\mu_0)$ settles down to a solution of the RDE. Note that as well as finding solutions of a given RDE, this method can be used to study endogeny via Theorem 11(c).

8.2. *Smoothness of densities.* For linear RDEs, under appropriate assumptions one can show that fixed points are unique and have $C^\infty$ densities, and use this as a basis for a theoretical "exact sampling" scheme; see [26]. In the Quicksort example (21), smoothness of densities has been studied in [31]. It would be interesting to seek general smoothness results for solutions of max-type RDE.

8.3. *Dependence on parameters.* When an RDE depends on a parameter (e.g., our (77) and (96)); see also examples involving multiplicative martingales for branching processes, e.g., [23], Theorem 3.3), it is natural to ask whether the solution depends continuously on the parameter. This has apparently not been studied in any generality.

8.4. *Continuous-time analogs.* We set up recursive tree processes as discrete-generation processes, analogous to discrete-time Markov chains. Let us mention two contexts where continuous-parameter analogs of RTPs arise. The first is the classical KPP equation, which is the analog of (54) for branching Brownian motion; see [37] for a recent probabilistic account. The second concerns the maximum $X$ of standard Brownian excursion of length 1. By scaling, the maximum $X_t$ for Brownian excursion of length $t$ satisfies $X_t \overset{d}{=} t^{1/2}X$. Since excursions above higher levels are independent (conditional on their lengths), we can write (for infinitesimal $\delta$)

$$X = \delta + \max_i t_i^{1/2}(\delta) \ X_i$$

where $(t_i(\delta), i \geq 1)$ are the lengths of excursions above level $\delta$ within standard Brownian excursion. See [18] for this kind of decomposition.

8.5. *Process-valued analogs.* There are examples where the distribution arising in an RDE is the distribution of a stochastic process, rather than a single real-valued random variable. Here is an illustration.

*Birth and assassination process* [9]. Start with one individual at time 0. During each individual's lifetime, children are born at the times of a Poisson (rate $\lambda$) process. An individual cannot die before the time of its parent's death (time 0, for the original individual); after that time, the individual lives for a further random time $S$, i.i.d. over individuals. Consider the random time $X$ at which the process becomes extinct. It is not hard to show [9] that $X < \infty$ a.s. under the assumption $\inf_{u>0} u^{-1}\mathbb{E}\exp(uS) < 1/\lambda$.

It does not seem possible to find an equation for $X$ itself, but one can study a process $(X(t), 0 \leq t < \infty)$ for which $X = X(0)$. Specifically, first set up the process of all possible descendants; for a realization, for each $t$ let $X(t)$ be time until extinction, in the modification where the first individual has a "fictional parent" who dies at time $t$. One can now argue that the process $(X(t))$ satisfies the RDE

$$X(t) \overset{d}{=} t + S + \max_{i \,:\, \xi_i \leq t + S}(\xi_i + X_i(t + S - \xi_i))$$

where $(\xi_i)$ are the points of a Poisson rate $\lambda$ process on $(0, \infty)$. This RDE has not been studied.

OPEN PROBLEM 65. Study the scaling behavior of $X$ in the limit as $1/\lambda \downarrow \inf_u u^{-1} \mathbb{E} \exp(uS)$.

8.6. *Matchings in random regular graphs.* Fix $r \geq 2$ and consider a random $r$-regular $n$-vertex graph $\mathcal{G}(n, r)$. Attach independent exponential(1) distributed random weights $(\xi_e)$ to edges. As in Section 4.6, let $M(n, r)$ be the maximum weight of a partial matching in $\mathcal{G}(n, r)$. The $n \to \infty$ limit of $\mathcal{G}(n, r)$, in the sense of local weak convergence, is the infinite $r$-regular tree $\mathbb{T}_r$. Thus one can seek to carry out the general program formalizing the cavity method (Section 7.5) in this setting. Recent work [33] provides interesting positive and negative results. The relevant RDE is [cf. (51)]

$$(99) \qquad X \stackrel{d}{=} \max_{1 \leq i \leq r-1} (0, \xi_i - X_i) \qquad (S = \mathbb{R}^+).$$

THEOREM 66 ([33]). *Let $T_{r-1}$ be the map associated with the RDE (99). Then $T_{r-1}^2$ has a unique invariant distribution. Moreover, for $(X_i)$ with the invariant distribution,*

$$\lim_n n^{-1} \mathbb{E} M(n, r)$$

$$= \frac{1}{2} \mathbb{E} \sum_{i=1}^r \xi_i \mathbb{1}\left(\xi_i - X_i = \max_{1 \leq j \leq r}(\xi_j - X_j) > 0\right)$$

$$= \frac{rb^{r-1}}{2} \int_0^\infty t e^{-t} (1 - e^{-t}(1-b))^{r-1} \, dt$$

$$\quad + \frac{r(r-1)(1-b)}{2}$$

$$\quad \times \int_0^\infty \int_0^t t e^{-t} e^{-z} (1 - e^{-z}(1-b))^{r-2} (1 - e^{-t+z}(1-b))^{r-1} \, dz \, dt$$

*where $b$ is the unique solution of $b = 1 - \frac{1-b^r}{r(1-b)}$.*

Similar results for matchings on the sparse random graph model are also derived in [33]. But in the "dual" problem for independent subsets the relevant RDE turns out to have nonunique solution for $r \geq 5$, and nonuniqueness holds also for independent sets in the sparse random graph model above a certain critical value. So this setting provides an important test bed for the range of applicability of the method.

8.7. *Random fractal graphs.* A recent thesis [42] studies RDEs arising in the context of constructing random fractal graphs, and discusses examples such as the following:

$$X \stackrel{d}{=} g(X_1, X_2, \xi) \qquad (S = \mathbb{R}^+)$$

where

$$\xi \stackrel{d}{=} \text{Bern}(p),$$

$$g(x_1, x_2, 0) = 2\min(x_1, x_2),$$

$$g(x_1, x_2, 1) = \tfrac{1}{2}\max(x_1, x_2).$$

However, the form of functions $g$ are chosen for mathematical convenience, rather than being derived from an underlying richer model as in our examples.

8.8. *List of open problems and conjectures.* These can be fitted into four categories.

*Weakening hypotheses in general theorems.*
Open Problem 12: bivariate uniqueness implies endogeny.
Open Problem 19: solving linear RDE on $\mathbb{R}$.
Open Problem 18: endogeny for linear RDE on $\mathbb{R}^+$.
Open Problem 31: finiteness of discounted tree-sums.
Open Problem 45: convergence of BRW extremes.
Open Problem 49: nonendogeny of extremes in BRW.

*Existence and uniqueness of solution of particular RDEs.* Here one can also ask about endogeny.
Open Problem 46: extremes of branching random walk.
Open Problem 63: mean-field matching, $d \neq 1$.
Open Problem 64: mean-field TSP.
All three RDEs in Section 7.6.

*Endogeny and nonendogeny.*
Conjecture 34: discounted tree-sums.
Open Problem 62: domain of attraction for minimum matching.

*Scaling exponents.*
Open Problem 30: range of BRW and speed of algorithmic BRW.
Open Problem 65: birth and assassination.
All three RDEs in Section 7.6.

9. **Conclusion.** Here we attempt to review the big picture.

1. RDEs in general, and max-type RDEs in particular, arise in the study of a wide range of underlying stochastic models. Look again at the list of models in Table 1.

2. While for linear RDEs one has hope of general theory, the diversity of forms of the function $g(\cdot)$ listed in Table 1 makes it hard to envisage a general theory which encompasses many max-type examples in one axiomatic framework. Indeed it is not clear how to make any useful classification of our examples; we have given only an informal "simple/not simple" classification (start of Sections 4 and 5) based on whether there is a relatively easy a.s. construction of random variables satisfying the RDE.

3. The cavity method (Section 7.5) provides a range of examples new to the probability community. Existence and uniqueness of solutions has been proved rigorously only in the special settings of mean-field matching in pseudo-dimension 1 (Section 7.3) and matchings and independent sets in random graphs (Section 8.6). It remains a challenge to establish existence and uniqueness in the other examples of Sections 7.4 and 7.6.

4. What is new in this survey is the discussion of endogeny, both the (fairly straightforward) general theory in Section 2, and the analysis of examples. In some cases one can study endogeny in a model parameterized by a distribution $\xi$ (e.g., Corollaries 17 and 26; Proposition 48) but in other cases (Theorems 21, 55 and 61) the argument relies on analytic calculations based on knowing an explicit formula for the invariant distribution for a specific $\xi$. For making the cavity method rigorous, one would like techniques to establish endogeny without knowing such explicit formulas.

DEPARTMENT OF STATISTICS
UNIVERSITY OF CALIFORNIA
367 EVANS HALL #3860
BERKELEY, CALIFORNIA 94720
USA
E-MAIL: aldous@stat.berkeley.edu

INSTITUTE FOR MATHEMATICS
AND ITS APPLICATIONS
UNIVERSITY OF MINNESOTA
400 LIND HALL
207 CHURCH STREET
MINNEAPOLIS, MINNESOTA 55414
USA
E-MAIL: antar@ima.umn.edu